\title{Spectral Measures and Generating Series for Nimrep Graphs in Subfactor Theory {II}: $SU(3)$}
\author{
        David E. Evans and Mathew Pugh \\ \\
        School of Mathematics, \\
        Cardiff University, \\
        Senghennydd Road, \\
        Cardiff, CF24 4AG, \\
        Wales, U.K.
}
\date{\today}

\documentclass[12pt]{article}
\usepackage{amssymb}
\usepackage{graphicx}

\newtheorem{Def}{Definition}[section]

\newtheorem{Thm}[Def]{Theorem}

\begin{document}
\maketitle

\begin{abstract}
We complete the computation of spectral measures for $SU(3)$ nimrep graphs arising in subfactor theory, namely the $SU(3)$ $\mathcal{ADE}$ graphs associated with $SU(3)$ modular invariants and the McKay graphs of finite subgroups of $SU(3)$.
For the $SU(2)$ graphs the spectral measures distill onto very special subsets of the semicircle/circle, whilst for the $SU(3)$ graphs the spectral measures distill onto very special subsets of the discoid/torus. The theory of nimreps allows us to compute these measures precisely. We have previously determined spectral measures for some nimrep graphs arising in subfactor theory, particularly those associated with all $SU(2)$ modular invariants, all subgroups of $SU(2)$, the torus $\mathbb{T}^2$, $SU(3)$, and some $SU(3)$ graphs.
\end{abstract}

\section{Introduction} \label{sect:intro}

The Verlinde algebra of $SU(n)$ at level $k$ is represented by a non-degenerately braided system of endomorphisms ${}_N \mathcal{X}_N$ on a type $\mathrm{III}_1$ factor $N$ with fusion rules $\lambda \mu = \bigoplus_{\nu} N_{\lambda \nu}^{\mu} \nu$ \cite{wassermann:1998}. The fusion matrices $N_{\lambda} = [N_{\rho \lambda}^{\sigma}]_{\rho,\sigma}$ are a family of commuting normal matrices, and themselves give a representation of the fusion rules of the positive energy representations of the loop group of $SU(n)$ at level $k$, $N_{\lambda} N_{\mu} = \sum_{\nu} N_{\lambda \nu}^{\mu} N_{\nu}$, the regular representation.
This family $\{ N_{\lambda} \}$ of fusion matrices can be simultaneously diagonalised:
\begin{equation} \label{eqn:verlinde_formula}
N_{\lambda} = \sum_{\sigma} \frac{S_{\sigma, \lambda}}{S_{\sigma,1}} S_{\sigma} S_{\sigma}^{\ast},
\end{equation}
where $1$ is the trivial representation, and the eigenvalues $S_{\sigma, \lambda}/S_{\sigma,1}$ and eigenvectors $S_{\sigma} = [S_{\sigma, \mu}]_{\mu}$ are described by the symmetric modular $S$ matrix.

Braided subfactors $N \subset M$ (the dual canonical endomorphism is in $\Sigma({}_N \mathcal{X}_N)$, i.e. decomposes as a finite linear combination of endomorphisms in ${}_N \mathcal{X}_N$) yield modular invariants through the procedure of $\alpha$-induction which allows two extensions of $\lambda$ on $N$ to endomorphisms $\alpha^{\pm}_{\lambda}$ of $M$, such that the matrix $Z_{\lambda,\mu} = \langle \alpha_{\lambda}^+, \alpha_{\mu}^- \rangle$ is a modular invariant \cite{bockenhauer/evans/kawahigashi:1999, bockenhauer/evans:2000, evans:2003}.
The action of the $N$-$N$ sectors ${}_N \mathcal{X}_N$ on the $M$-$N$ sectors ${}_M \mathcal{X}_N$ produces a nimrep (non-negative matrix integer representation of the fusion rules)
$$G_{\lambda} G_{\mu} = \sum_{\nu} N_{\lambda \nu}^{\mu} G_{\nu}$$
whose spectrum reproduces exactly the diagonal part of the modular invariant, i.e.
\begin{equation} \label{eqn:verlinde_formulaG}
G_{\lambda} = \sum_i \frac{S_{i,\lambda}}{S_{i,1}} \psi_i \psi_i^{\ast},
\end{equation}
with the spectrum of $G_{\lambda} = \{ S_{\mu, \lambda}/S_{\mu,1}$ with multiplicity $Z_{\mu,\mu} \}$ \cite{bockenhauer/evans/kawahigashi:2000}. The labels $\mu$ of the non-zero diagonal elements are called the exponents of $Z$, counting multiplicity.

Every $SU(2)$ and $SU(3)$ modular invariant can be realised by $\alpha$-induction for a suitable braided subfactor \cite{ocneanu:2000ii, ocneanu:2002, xu:1998, bockenhauer/evans:1999i, bockenhauer/evans:1999ii, bockenhauer/evans/kawahigashi:1999, bockenhauer/evans/kawahigashi:2000}, \cite{ocneanu:2000ii, ocneanu:2002, evans/pugh:2009ii} respectively. For $SU(2)$, the classification of Cappelli, Itzykson and Zuber \cite{cappelli/itzykson/zuber:1987ii} of $SU(2)$ modular invariants is understood in the following way. Suppose $N \subset M$ is a braided subfactor which realises the modular invariant $Z_{\mathcal{G}}$. Evaluating the nimrep $G$ at the fundamental representation $\rho$, we obtain for the inclusion $N \subset M$ a matrix $G_{\rho}$, which is the adjacency matrix for the $ADE$ graph $\mathcal{G}$ which labels the modular invariant.
Since these $ADE$ graphs can be matched to the affine Dynkin diagrams -- the McKay graphs of the finite subgroups of $SU(2)$ -- di Francesco and Zuber \cite{di_francesco/zuber:1990} were guided to find candidates for classifying graphs for $SU(3)$ modular invariants by first considering the McKay graphs of the finite subgroups of $SU(3)$ to produce a candidate list of $\mathcal{ADE}$ graphs whose spectra described the diagonal part of the modular invariant. The classification of $SU(3)$ modular invariants was shown to be complete by Gannon \cite{gannon:1994}, and the complete list is given in \cite{evans/pugh:2009ii}. Ocneanu claimed \cite{ocneanu:2000ii, ocneanu:2002} that all $SU(3)$ modular invariants were realised by subfactors and this was shown in \cite{evans/pugh:2009ii}.
The figures for the list of the $\mathcal{ADE}$ graphs are given in \cite{behrend/pearce/petkova/zuber:2000}, or in \cite{evans/pugh:2009i, evans/pugh:2009ii}. However this list of nimreps has not been shown to be complete.

In general, different inclusions which yield different nimreps may still realise the same modular invariant, as is the case in $SU(3)$ with the inclusions for the graphs $\mathcal{E}_1^{(12)}$ and $\mathcal{E}_2^{(12)}$, which both realise the modular invariant $Z_{\mathcal{E}^{(12)}}$
\cite[Section 8]{bockenhauer/evans:2002}. Thus any modular invariant may have more than one nimrep associated to it (although this is not the case in $SU(2)$). However, in $SU(3)$ there is uniqueness in the reverse direction, that is, each nimrep has an unique modular invariant associated to it, due to the coincidence that at any level $k$ each $SU(3)$ modular invariant has a different trace.
Unlike the situation for $SU(2)$, there is a mismatch between the list of nimreps associated to each $SU(3)$ modular invariant and the McKay graphs of the finite subgroups of $SU(3)$ which are also the nimreps of the representation theory of the group. The latter also have a diagonalisation as in (\ref{eqn:verlinde_formula}), with diagonalising matrix $S = \{ S_{ij} \}$ usually non-symmetric, where $i$ labels conjugacy classes and $j$ the irreducible characters (see \cite[Section 8.7]{evans/kawahigashi:1998} and \cite[Section 4]{evans/pugh:2009v}).  Both of these kinds of nimreps will play a role in this paper.

In \cite{evans/pugh:2009v} we determined spectral measures for some nimrep graphs arising in subfactor theory, particularly those associated with all $SU(2)$ modular invariants and all subgroups of $SU(2)$. Our methods gave an alternative approach to deriving the results of Banica and Bisch \cite{banica/bisch:2007} for $ADE$ graphs and subgroups of $SU(2)$, and explained the connection between their results for affine $ADE$ graphs and the Kostant polynomials. We also determined spectral measures for the torus $\mathbb{T}^2$ and $SU(3)$, and some $SU(3)$ graphs, namely $\mathcal{A}^{(n)}$, $\mathcal{D}^{(3k)}$ and $\mathcal{A}^{(n)\ast}$, for integers $n \geq 4$, $k \geq 2$. We now complete the computation of the spectral measures for the $SU(3)$ $\mathcal{ADE}$ graphs in this present work, as well as all finite subgroups of $SU(3)$.

Suppose $A$ is a unital $C^{\ast}$-algebra with state $\varphi$.
If $b \in A$ is a normal operator then there exists a compactly supported probability measure $\mu_b$ on the spectrum $\sigma(b) \subset \mathbb{C}$ of $b$, uniquely determined by its moments
\begin{equation} \label{eqn:moments_normal_operator}
\varphi(b^m b^{\ast n}) = \int_{\sigma(b)} z^m \overline{z}^n \mathrm{d}\mu_b (z),
\end{equation}
for non-negative integers $m$, $n$.

We computed in \cite{evans/pugh:2009v} such spectral measures and generating series when $b$ is the normal operator $\Delta = G_{\rho}$ acting on the Hilbert space of square summable functions on the graph, for the nimreps $G$ described above, i.e. $G_{\rho}$ is the adjacency matrix of the $ADE$ and affine $ADE$ graphs in $SU(2)$ and certain graphs in $SU(3)$. We computed the spectral measure for the vacuum, i.e. the distinguished vertex of the graph which has lowest Perron-Frobenius weight. However the spectral measures for the other vertices of the graph could also be computed by the same methods.

In particular, for $SU(2)$, we can understand the spectral measures for the torus $T$ and $SU(2)$ as follows.
If $w_Z$ and $w_N$ are the self adjoint operators arising from the McKay graph of the fusion rules of the representation theory of $T$ and $SU(2)$, then the spectral measures in the vacuum state can be describe in terms of semicircular law, on the interval $[-2,2]$ which is the spectrum of either
as the image of the map $z \in T \rightarrow z + z^{-1}$ \cite[Sections 2 \& 3.1]{evans/pugh:2009v}:
$$\mathrm{dim}\left( \left(\otimes^k M_2 \right)^{\mathbb{T}} \right) \;\; = \;\; \varphi(w_Z^{2k}) \;\; = \;\; \frac{1}{\pi} \int_{-2}^2 x^{2k} \frac{1}{\sqrt{4-x^2}} \; \mathrm{d}x \, ,$$
$$\mathrm{dim}\left( \left(\otimes^k M_2 \right)^{SU(2)} \right) \;\; = \;\; \varphi(w_N^{2k}) \;\; = \;\; \frac{1}{2\pi} \int_{-2}^2 x^{2k} \sqrt{4-x^2} \; \mathrm{d}x \, .$$
The fusion matrix for $SO(3)$ is just $\mathbf{1} + \Delta$, where $\Delta$ is the fusion matrix for $SU(2)$, and thus is equal to the infinite $SU(3)$ $\mathcal{ADE}$ graph $\mathcal{A}^{(\infty)\ast}$. Thus the spectral measure $\mu$ in the vacuum state (over $[-1,3]$) for $SO(3)$ has semicircle distribution with mean 1 and variance 1, i.e. $\mathrm{d}\mu(x) = \sqrt{4 - (x-1)^2} \mathrm{d}x$ \cite[Section 7.3]{evans/pugh:2009v}.

The spectral weight for $SU(2)$ arises from the Jacobian of a change of variable between the interval $[-2,2]$ and the circle.
Then for $T^2$ and $SU(3)$, the 3-cusp discoid $\mathfrak{D}$ in the complex plane is the image of the two-torus under the map $(\omega_1, \omega_2) \rightarrow \omega_1 + \omega_2^{-1} + \omega_1^{-1} \omega_2$, which is the spectrum of the corresponding normal operators on the Hilbert spaces of the fusion graphs. The corresponding spectral measures are then described by a corresponding Jacobian or discriminant as \cite[Theorems 3 \& 5]{evans/pugh:2009v}:
$$\mathrm{dim}\left( \left(\otimes^k M_3 \right)^{\mathbb{T}^2} \right) \;\; = \;\; \varphi(|v_Z|^{2k}) \;\; = \;\; \frac{3}{\pi^2} \int_{\mathfrak{D}} |z|^{2k} \frac{1}{\sqrt{27 - 18z\overline{z} + 4z^3 + 4\overline{z}^3 -z^2 \overline{z}^2}} \; \mathrm{d}z \, ,$$
$$\mathrm{dim}\left( \left(\otimes^k M_3 \right)^{SU(3)} \right) \;\; = \;\; \varphi(|v_N|^{2k}) \;\; = \;\; \frac{1}{2\pi^2} \int_{\mathfrak{D}} |z|^{2k} \sqrt{27 - 18z\overline{z} + 4z^3 + 4\overline{z}^3 -z^2 \overline{z}^2} \; \mathrm{d}z \, ,$$
where $\mathrm{d}z:=\mathrm{d}\,\mathrm{Re}z \; \mathrm{d}\,\mathrm{Im}z$ denotes the Lebesgue measure on $\mathbb{C}$.

For the $SU(2)$ and $SU(3)$ graphs, the spectral measures distill onto very special subsets of the semicircle/circle ($SU(2)$) and discoid/torus ($SU(3)$), and the theory of nimreps allows us to compute these measures precisely. In the present work we complete the computation of the spectral measures for the $SU(3)$ $\mathcal{ADE}$ graphs in Section \ref{sect:spec_measureSU(3)ADE}, and compute the spectral measures for the finite subgroups of $SU(3)$ in Section \ref{sect:spec_measure-subgroupsSU(3)}.

\section{Spectral Measures of the $SU(3)$ $\mathcal{ADE}$ Graphs} \label{sect:spec_measureSU(3)ADE}

Let $\Delta_{\mathcal{G}}$ denote the adjacency matrix of a finite graph $\mathcal{G}$ for which $\Delta_{\mathcal{G}}$ is normal, and let $e_v$ denote the basis vector in $\ell^2(\mathcal{G})$ corresponding to the vertex $v$ of $\mathcal{G}$. The inner product $\langle \Delta_{\mathcal{G}}^m (\Delta_{\mathcal{G}}^{\ast})^n e_v, e_v \rangle$ defines a spectral measure $\mu$ of $(\mathcal{G},v)$ which has $m,n^{\mathrm{th}}$ moment $\int z^m \overline{z}^n \mathrm{d}\mu(z) = \langle \Delta_{\mathcal{G}}^m (\Delta_{\mathcal{G}}^{\ast})^n e_v, e_v \rangle$. In this work we will compute the spectral measure in the vacuum state, i.e. when $v$ is the distinguished vertex $\ast$ of $\mathcal{G}$ which has lowest Perron-Frobenius weight. However the same method will work for any vertex $v$ of $\mathcal{G}$.
For convenience we will use the notation
\begin{equation} \label{def:Rmn}
R_{m,n}(\omega_1,\omega_2) := (\omega_1 + \omega_2^{-1} + \omega_1^{-1}\omega_2)^m (\omega_1^{-1} + \omega_2 + \omega_1\omega_2^{-1})^n,
\end{equation}
so that $\int_{\mathbb{T}^2} R_{m,n}(\omega_1,\omega_2) \mathrm{d}\varepsilon(\omega_1,\omega_2) = \int z^m \overline{z}^n \mathrm{d}\mu(z) = \langle \Delta_{\mathcal{G}}^m (\Delta_{\mathcal{G}}^{\ast})^n e_1, e_1 \rangle$.

Let $\beta^j$ be the eigenvalues of $\mathcal{G}$, with corresponding eigenvectors $x^j$, $j=1,\ldots,s$, where $s$ is the number of vertices of $\mathcal{G}$. Then as for $SU(2)$ \cite[Section 3]{evans/pugh:2009v}, $\Delta_{\mathcal{G}}^m (\Delta_{\mathcal{G}}^{\ast})^n = \mathcal{U} \Lambda_{\mathcal{G}}^m (\Lambda_{\mathcal{G}}^{\ast})^n \mathcal{U}^{\ast}$, where $\Lambda_{\mathcal{G}}$ is the diagonal matrix $\Lambda_{\mathcal{G}} = \mathrm{diag}(\beta^1, \beta^2, \ldots, \beta^s)$ and $\mathcal{U} = (x^1, x^2, \ldots, x^s)$, so that
\begin{eqnarray}
\int_{\mathbb{T}^2} R_{m,n}(\omega_1,\omega_2) \mathrm{d}\varepsilon(\omega_1,\omega_2) & = & \langle \mathcal{U} \Lambda_{\mathcal{G}}^m (\Lambda_{\mathcal{G}}^{\ast})^n \mathcal{U}^{\ast} e_1, e_1 \rangle \;\; = \;\; \langle \Lambda_{\mathcal{G}}^m (\Lambda_{\mathcal{G}}^{\ast})^n \mathcal{U}^{\ast} e_1, \mathcal{U}^{\ast} e_1 \rangle \nonumber \\
& = & \sum_{j=1}^s (\beta^j)^m (\overline{\beta^j})^n |y_j|^2, \label{eqn:moments_general_SU(3)graph}
\end{eqnarray}
where $y_j = x^j_1$ is the first entry of the eigenvector $x^j$.

\begin{figure}[tb]
\begin{center}
  \includegraphics[width=70mm]{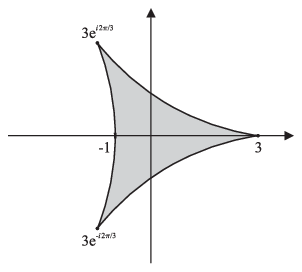}\\
  \caption{The (3-cusp) discoid $\mathfrak{D}$, the union of the deltoid and its interior.} \label{fig:hypocycloid-S}
\end{center}
\end{figure}

Now suppose $\mathcal{G}$ is a finite $\mathcal{ADE}$ graph with distinguished vertex $\ast$ which is the vertex with lowest Perron-Frobenius weight. Every eigenvalue $\beta^{(\lambda)}$ of $\mathcal{G}$ is a ratio of the $S$-matrix given by $\beta^{(\lambda)} = S_{\rho \lambda}/S_{0 \lambda}$, for a Dynkin label $\lambda$, with corresponding eigenvector $(\psi^{\lambda}_a)_{a \in \mathfrak{V}(\mathcal{G})}$.
Suppose $\mathcal{G}$ is the nimrep given by a braided subfactor which realises the modular invariant $Z$. Then the spectrum of the adjacency matrix $\Delta_{\mathcal{G}}$ of $\mathcal{G}$ is $\sigma(\Delta_{\mathcal{G}}) = \{ \beta^{(\lambda)} | \; \lambda \in \mathrm{Exp} \}$, where $\mathrm{Exp}$ is the set of exponents of $Z$ counting multiplicity.
The moments $\int_{\mathbb{T}^2} R_{m,n}(\omega_1,\omega_2) \mathrm{d}\varepsilon(\omega_1,\omega_2)$ of $\mathcal{G}$ over $\mathbb{T}^2$ are given by (\ref{eqn:moments_general_SU(3)graph}) or \cite[Equation (48)]{evans/pugh:2009v}:
\begin{equation} \label{eqn:moments_SU(3)}
\int_{\mathbb{T}^2} R_{m,n}(\omega_1,\omega_2) \mathrm{d}\varepsilon(\omega_1,\omega_2) = \sum_{\lambda \in \mathrm{Exp}} (\beta^{(\lambda)})^m (\overline{\beta^{(\lambda)}})^n |\psi^{\lambda}_{\ast}|^2.
\end{equation}
The spectrum $\sigma(\Delta_{\mathcal{G}})$ of the adjacency matrix $\Delta_{\mathcal{G}}$ of any $SU(3)$ $\mathcal{ADE}$ graph $\mathcal{G}$ is contained in the spectrum $\sigma(\Delta) = \mathfrak{D} = \{ \omega_1 + \omega_2^{-1} + \omega_1^{-1}\omega_2 | \; \omega_1,\omega_2 \in \mathbb{T} \}$ of the adjacency matrix $\Delta$ of $\mathcal{A}^{(\infty)}$. The (3-cusp) discoid $\mathfrak{D}$ is the surface given by the union of the deltoid and its interior, illustrated in Figure \ref{fig:hypocycloid-S}.

Thus the support of the probability measure $\mu$ of $\mathcal{G}$ is contained in the discoid $\mathfrak{D}$.
There is a map $\Phi:\mathbb{T}^2 \rightarrow \mathfrak{D}$ from the torus onto $\mathfrak{D}$ given by
\begin{equation} \label{Phi:T^2->S}
\Phi(\omega_1,\omega_2) = \omega_1 + \omega_2^{-1} + \omega_1^{-1} \omega_2,
\end{equation}
where $\omega_1, \omega_2 \in \mathbb{T}$.

The Weyl group of $SU(3)$ is the permutation group $S_3$. Consider the group $S_3$ as a subgroup of $GL(2,\mathbb{Z})$, generated by the matrices $T_2$, $T_3$, of orders 2, 3 respectively, given by
\begin{equation} \label{T1,T2}
T_2 = \left( \begin{array}{cc} 0 & -1 \\ -1 & 0 \end{array} \right), \qquad T_3 = \left( \begin{array}{cc} 0 & -1 \\ 1 & -1 \end{array} \right).
\end{equation}
The action of $S_3$ on $\mathbb{T}^2$ is given by $T(\omega_1,\omega_2) = (\omega_1^{a_{11}} \omega_2^{a_{12}}, \omega_1^{a_{21}} \omega_2^{a_{22}})$, for $T = (a_{ij}) \in S_3$. For $(\omega_1,\omega_2) = (e^{2 \pi i \theta_1}, e^{2 \pi i \theta_2})$, we will define the action of $S_3$ on $(\theta_1, \theta_2) \in [0,1] \times [0,1]$ by $T(\theta_1, \theta_2) = (a_{11} \theta_1 + a_{12} \theta_2, a_{21} \theta_1 + a_{22} \theta_2)$ for $T = (a_{ij}) \in S_3$ (notice that $T(e^{2 \pi i \theta_1}, e^{2 \pi i \theta_2}) = (e^{2 \pi i \theta_1'}, e^{2 \pi i \theta_2'})$ where $(\theta_1', \theta_2') = T(\theta_1, \theta_2)$).

\begin{figure}[tb]
\begin{center}
  \includegraphics[width=55mm]{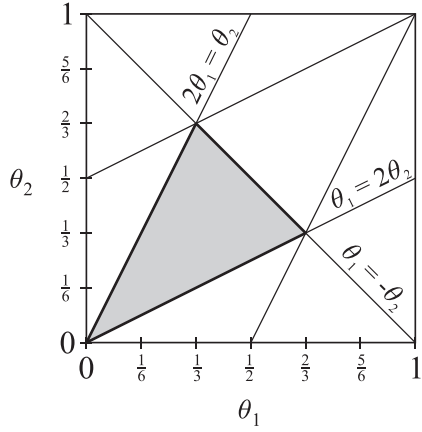}\\
 \caption{A fundamental domain $C$ of $\mathbb{T}^2/S_3$.} \label{fig:poly-4}
\end{center}
\end{figure}

The quotient $\mathbb{T}^2/S_3$ is topologically homeomorphic to the discoid $\mathfrak{D}$ \cite[Section 6.1]{evans/pugh:2009v}. The deltoid, which is the boundary of the discoid $\mathfrak{D}$, is given by the lines $\theta_1 = 1 - \theta_2$, $\theta_1 = 2 \theta_2$ and $2 \theta_1 = \theta_2$. The diagonal $\theta_1 = \theta_2$ in $\mathbb{T}^2$ is mapped to the real interval $[-1,3] \subset \mathfrak{D}$. A fundamental domain $C$ of $\mathbb{T}^2$ under the action of the group $S_3$ is illustrated in Figure \ref{fig:poly-4}, where the axes are labelled by the parameters $\theta_1$, $\theta_2$ in $(e^{2 \pi i \theta_1}, e^{2 \pi i \theta_2}) \in \mathbb{T}^2$. The boundaries of $C$ map to the deltoid. The torus $\mathbb{T}^2$ contains six copies of $C$.

Any probability measure $\varepsilon$ on $\mathbb{T}^2$ produces a probability measure $\mu = \Phi^{\ast}\varepsilon$ on $\mathfrak{D}$. There is a bijection between $S_3$-invariant probability measures on $\mathbb{T}^2$ and probability measures on $\mathfrak{D}$.
We will compute $S_3$-invariant spectral measures of the $SU(3)$ graphs on $\mathbb{T}^2$. It was shown in \cite[Section 7.1]{evans/pugh:2009v} that the eigenvalues $\beta^{(\lambda)} \in \mathfrak{D}$, $\lambda \in \mathrm{Exp}$, of an $SU(3)$ $\mathcal{ADE}$ graph $\mathcal{G}$ are given by
\begin{equation} \label{eqn:theta->lambda}
\beta^{(\lambda)} = \Phi(e^{2\pi i \theta_1},e^{2\pi i \theta_2}),
\end{equation}
where $\theta_1 = (\lambda_1 + 2\lambda_2 + 3)/3n$ and $\theta_2 = (2\lambda_1 + \lambda_2 + 3)/3n$.

Under the change of variable $z = e^{2 \pi i \theta_1} + e^{-2 \pi i \theta_2} + e^{2 \pi i (\theta_2 - \theta_1)}$, we have
\begin{eqnarray*}
x & := & \mathrm{Re}(z) \;\; = \;\; \cos(2 \pi \theta_1) + \cos(2 \pi \theta_2) + \cos(2 \pi (\theta_2 - \theta_1)), \\
y & := & \mathrm{Im}(z) \;\; = \;\; \sin(2 \pi \theta_1) - \sin(2 \pi \theta_2) + \sin(2 \pi (\theta_2 - \theta_1)).
\end{eqnarray*}
Then
\begin{eqnarray}
\lefteqn{ \int_{\mathbb{T}^2} (\omega_1 + \omega_2^{-1} + \omega_1^{-1}\omega_2)^m (\omega_1^{-1} + \omega_2 + \omega_1\omega_2^{-1})^n \mathrm{d}\omega_1 \; \mathrm{d}\omega_2 } \nonumber \\
& = & 6 \int_{\mathfrak{D}} (x+iy)^m (x+iy)^n |J|^{-1} \mathrm{d}x \; \mathrm{d}y, \label{eqn:integral-A(infty)6_over_S}
\end{eqnarray}
where the Jacobian $J = \mathrm{det}(\partial(x,y)/\partial(\theta_1,\theta_2))$ is the determinant of the Jacobian matrix.
The Jacobian $J = J (\theta_1,\theta_2)$ is given by \cite[Equation (39)]{evans/pugh:2009v}
\begin{equation} \label{eqn:J[theta]}
J (\theta_1,\theta_2) = 4 \pi^2 (\sin(2 \pi (\theta_1 + \theta_2)) - \sin(2 \pi (2\theta_1 - \theta_2)) - \sin(2 \pi (2\theta_2 - \theta_1))).
\end{equation}
The Jacobian is real and vanishes on the deltoid, the boundary of the discoid $\mathfrak{D}$. For the values of $\theta_1$, $\theta_2$ such that $(e^{2 \pi i \theta_1},e^{2 \pi i \theta_2})$ are in the interior of the fundamental domain $C$ illustrated in Figure \ref{fig:poly-4}, the value of $J$ is always negative. In fact, restricting to any one of the fundamental domains shown in Figure \ref{fig:poly-4}, the sign of $J$ is constant. It is negative over three of the fundamental domains, and positive over the remaining three.
When evaluating $J$ at a point in $z \in \mathfrak{D}$, we pull back $z$ to $\mathbb{T}^2$. However, there are six possibilities for $(\omega_1,\omega_2) \in \mathbb{T}^2$ such that $\Phi(\omega_1,\omega_2) = z$, one in each of the fundamental domains of $\mathbb{T}^2$ in Figure \ref{fig:poly-4}. Thus over $\mathfrak{D}$, $J$ is only determined up to a sign.
To obtain a positive measure over $\mathfrak{D}$ we take the absolute value $|J|$ of the Jacobian in the integral (\ref{eqn:integral-A(infty)6_over_S}).

Since $J^2$ is invariant under the action of $S_3$ it can be written in terms of $z$, $\overline{z}$, namely $J (z,\overline{z})^2 = 4 \pi^4 (27 - 18z\overline{z} + 4z^3 + 4\overline{z}^3 -z^2 \overline{z}^2)$ for $z \in \mathfrak{D}$ \cite[Section 6.1]{evans/pugh:2009v}. Since $J$ is real, $J^2 \geq 0$, so we can write
\begin{equation} \label{eqn:J(z)}
|J (z,\overline{z})| = 2 \pi^2 \sqrt{27 - 18z\overline{z} + 4z^3 + 4\overline{z}^3 -z^2 \overline{z}^2}.
\end{equation}

\begin{figure}[tb]
\begin{center}
  \includegraphics[width=55mm]{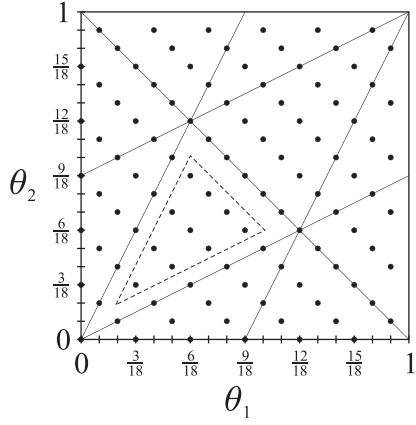}\\
 \caption{The points $(\theta_1,\theta_2)$ such that $(e^{2 \pi i \theta_1}, e^{2 \pi i \theta_2}) \in D_6$.} \label{fig:D6}
\end{center}
\end{figure}

In \cite[Theorem 4]{evans/pugh:2009v} it was shown that the spectral measure (on $\mathbb{T}^2$) for the graph $\mathcal{A}^{(n)}$ is the measure $J^2 \mathrm{d}^{(n)}$ (up to a factor of $16 \pi^4$), where $\mathrm{d}^{(n)}$ is the uniform measure on
\begin{equation} \label{def:Dl}
D_n = \{ (e^{2 \pi i q_1/3n}, e^{2 \pi i q_2/3n}) \in \mathbb{T}^2 | \; q_1,q_2 = 0, 1, \ldots, 3n-1; q_1 + q_2 \equiv 0 \textrm{ mod } 3 \}.
\end{equation}
The points $(\theta_1,\theta_2) \in [0,1]^2$ for which $(e^{2 \pi i \theta_1}, e^{2 \pi i \theta_2}) \in D_6$ are illustrated in Figure \ref{fig:D6}. Notice that the points in the interior of the fundamental domain $C$ (those enclosed by the dashed line) correspond to the vertices of the graph $\mathcal{A}^{(6)}$.

The Jacobian $J$ (\ref{eqn:J(z)}) will appear in Section \ref{sect:spec_measure-subgroupsSU(3)} as the discriminant in solutions for the inverse image $\Phi^{-1}(z) \in \mathbb{T}^2$ of $z \in \mathfrak{D}$.  This discriminant also appears in the work of Gepner \cite[Equation (2.64)]{gepner:1991} as the measure required to make the polynomials $S_{\mu}(z,\overline{z})$ orthogonal, where the polynomials $S_{\nu}(x,y)$ are defined by $S_{(0,0)}(x,y) = 1$, and $x S_{\nu}(x,y) = \sum_{\mu} \Delta_{\mathcal{A}}(\nu, \mu) S_{\mu}(x,y)$ and $y S_{\nu}(x,y) = \sum_{\mu} \Delta_{\mathcal{A}}^T(\nu, \mu) S_{\mu}(x,y)$ for vertices $\nu$ of $\mathcal{A}^{(n)}$.
The Jacobian may be written in terms of $(\omega_1, \omega_2) \in \mathbb{T}^2$ as \cite[Equation (40)]{evans/pugh:2009v}:
\begin{equation} \label{eqn:J[omega]}
i J (\omega_1,\omega_2)/2 \pi^2 = \omega_1 \omega_2 - \omega_1^{-1} \omega_2^{-1} - \omega_1^2 \omega_2^{-1} + \omega_1^{-2} \omega_2 - \omega_1^{-1} \omega_2^2 + \omega_1 \omega_2^{-2}.
\end{equation}

\noindent \emph{Remark:}
Let $1$ denote the trivial representation, $\rho$ the fundamental representation of $SU(3)$ and $\overline{\rho}$ its conjugate representation. Kuperberg \cite[Conjecture 3.4]{kuperberg:1994} conjectured that certain $A_2$-$(k,n)$-tangles in the sense of \cite[Section 2.3]{evans/pugh:2009iii} which do not contain elliptic faces in the sense of \cite[Section 4]{kuperberg:1996} are a basis for $\mathrm{Hom}(1,\rho^n \overline{\rho}^k)\,,$ and observed that  it is sufficient to show that the number of such $A_2$-$(k,n)$-tangles is given by the coefficient of the term $\omega_1^{-1}\omega_2^2$ in the polynomial:
\begin{equation} \label{eqn:KupA2}
(\omega_1 + \omega_2^{-1} + \omega_1^{-1}\omega_2)^k (\omega_1^{-1} + \omega_2 + \omega_1\omega_2^{-1})^n (\omega_1^{-1}\omega_2^2 - \omega_1\omega_2 + \omega_1^2\omega_2^{-1} - \omega_1\omega_2^{-2} + \omega_1^{-1}\omega_2^{-1} - \omega_1^{-2}\omega_2).
\end{equation}
In the special case of $k=n$, this follows from   \cite[Theorem 5]{evans/pugh:2009v}.
Denote by $p(k,n)$ the polynomial in (\ref{eqn:KupA2}). The coefficient $c$ of the term $\omega_1^{-1}\omega_2^2$ in $p(k,k)$ is given by the integral $\int_{\mathbb{T}^2} p(k,k) \omega_1\omega_2^{-2} \, \mathrm{d}\omega_1 \, \mathrm{d}\omega_2$. Averaging over the orbit of $\omega_1\omega_2^{-2}$ under the action of the Weyl group $S_3$ of $SU(3)$ gives $i J/2 \pi^2$ as in (\ref{eqn:J[omega]}), thus $c = \int_{\mathbb{T}^2} (\omega_1 + \omega_2^{-1} + \omega_1^{-1}\omega_2)^k (\omega_1^{-1} + \omega_2 + \omega_1\omega_2^{-1})^k J^2 \, \mathrm{d}\omega_1 \, \mathrm{d}\omega_2 / 24\pi^4$.
This is  the dimension of the path algebra $(\bigotimes^k M_3)^{SU(3)}$ \cite[Corollary 1]{evans/pugh:2009v}, which has basis given by the $A_2$-$(k,k)$-tangles which do not contain elliptic faces \cite[Lemma 2.12]{evans/pugh:2009iii}.

It was shown in \cite[Sections 7.4 \& 7.5]{evans/pugh:2009v} that the spectral measure for the graphs $\mathcal{E}^{(8)}$ and $\mathcal{E}_1^{(12)}$ cannot be written as a linear combination of measures of the form $\mathrm{d}^{(p)}$, $J^2 \mathrm{d}^{(p)}$, $J^2 \mathrm{d}_{p/2} \times \mathrm{d}_{p/2}$ and $\mathrm{d}_{p/2} \times \mathrm{d}_{p/2}$ for $p \in \mathbb{N}$, where $\mathrm{d}_p$ is the uniform measure on the $2p^{\mathrm{th}}$ roots of unity. However, if we introduce two new measures $\mathrm{d}^{((n))}$, $\mathrm{d}^{(n,k)}$, we can write the spectral measures for $\mathcal{E}^{(8)}$, $\mathcal{E}_1^{(12)}$ and the other $SU(3)$ $\mathcal{ADE}$ nimrep graphs as linear combinations of these measures.

\begin{Def} \label{def:4measures}
Let $\omega = e^{2 \pi i/3}$, $\tau = e^{2 \pi i/n}$. We define the following measures on $\mathbb{T}^2$:
\begin{itemize}
\item[(1)] $\mathrm{d}_m \times \mathrm{d}_n$, where $\mathrm{d}_k$ is the uniform measure on the $2k^{\mathrm{th}}$ roots of unity, for $k \in \mathbb{N}$.
\item[(2)] $\mathrm{d}^{(n)}$, the uniform measure on $D_n$ for $n \in \mathbb{N}$.
\item[(3)] $\mathrm{d}^{((n))}$, the uniform measure on the $S_3$-orbit of the points $(\tau, \tau)$, $(\overline{\omega} \, \overline{\tau}, \omega)$, $(\omega, \overline{\omega} \, \overline{\tau})$, for $n \in \mathbb{Q}$, $n \geq 2$.
\item[(4)] $\mathrm{d}^{(n,k)}$, the uniform measure on the $S_3$-orbit of the points $(\tau \, e^{2 \pi i k}, \tau)$, $(\tau, \tau \, e^{2 \pi i k})$, $(\overline{\omega} \, \overline{\tau}, \omega \, e^{2 \pi i k})$, $(\omega \, e^{2 \pi i k}, \overline{\omega} \, \overline{\tau})$, $(\overline{\omega} \, \overline{\tau} \, e^{-2 \pi i k}, \omega \, e^{-2 \pi i k})$, $(\omega \, e^{-2 \pi i k}, \overline{\omega} \, \overline{\tau} \, e^{-2 \pi i k})$, for $n,k \in \mathbb{Q}$, $n > 2$, $0 \leq k \leq 1/n$.
\end{itemize}
\end{Def}

\begin{figure}[tb]
\begin{minipage}[t]{7.5cm}
\begin{center}
  \includegraphics[width=55mm]{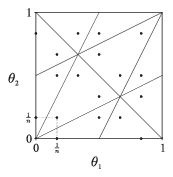}\\
 \caption{$Supp(\mathrm{d}^{((n))})$} \label{fig:poly-15}
\end{center}
\end{minipage}
\hfill
\begin{minipage}[t]{7.5cm}
\begin{center}
  \includegraphics[width=55mm]{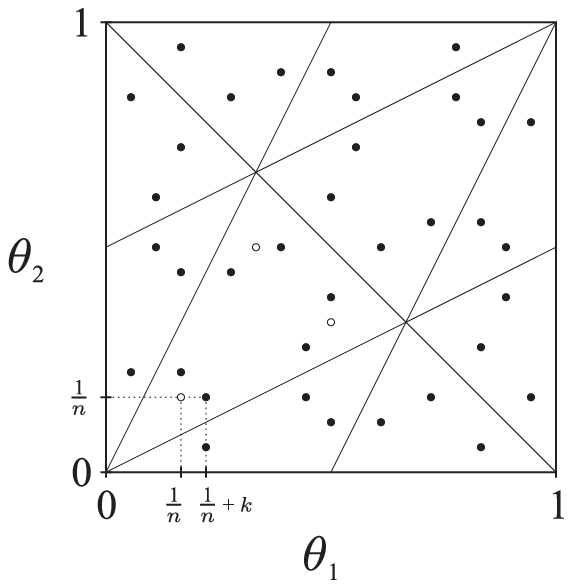}\\
 \caption{$Supp(\mathrm{d}^{(n,k)})$} \label{fig:poly-16}
\end{center}
\end{minipage}
\end{figure}

Let $Supp(\mathrm{d}\mu)$ denote the set of points $(\theta_1,\theta_2) \in [0,1]^2$ such that $(e^{2 \pi i \theta_1}, e^{2 \pi i \theta_2})$ is in the support of the measure $\mathrm{d}\mu$. The sets $Supp(\mathrm{d}^{((n))})$, $Supp(\mathrm{d}^{(n,k)})$ are illustrated in Figures \ref{fig:poly-15}, \ref{fig:poly-16} respectively. The white circles in Figure \ref{fig:poly-16} denote the points given by the measure $\mathrm{d}^{((n))}$. The cardinality $|Supp(\mathrm{d}_m \times \mathrm{d}_n)|$ of $Supp(\mathrm{d}_m \times \mathrm{d}_n)$ is $mn$, whilst $|Supp(\mathrm{d}^{(n)})| = |D_n| = 3n^2$ was shown in \cite[Section 7.1]{evans/pugh:2009v}. For $n > 2$ and $0 < k < 1/n$, $|Supp(\mathrm{d}^{((n))})| = 18$, whilst $|Supp(\mathrm{d}^{(n,k)})| = 36$. The cardinalities of the other sets are $|Supp(\mathrm{d}^{(n,0)})| = |Supp(\mathrm{d}^{(n,1/n)})| = 18$ for $n > 2$, and $|Supp(\mathrm{d}^{((2))})| = 9$.

It is easy to check that the following relations hold:
\begin{eqnarray}
3J^2 \; \mathrm{d}^{(3)} & = & J^2 \; \mathrm{d}_{3/2} \times \mathrm{d}_{3/2}, \label{eqn:measure-relation1} \\
\mathrm{d}^{((4))} & = & \frac{1}{24\pi^4} J^2 \; \mathrm{d}^{(4)}, \label{eqn:measure-relation2} \\
\mathrm{d}^{(n,0)} & = & \mathrm{d}^{((n))}, \label{eqn:measure-relation3} \\
\mathrm{d}^{(6,1/6)} & = & \mathrm{d}^{((2))} \;\; = \;\; \frac{1}{3}(4\mathrm{d}^{(2)} - \mathrm{d}^{(1)}). \label{eqn:measure-relation4}
\end{eqnarray}

\renewcommand{\arraystretch}{1.5}

\subsection{Graphs $\mathcal{D}^{(n)\ast}$}

Let $A$ denote the automorphism of order 3 on the exponents $\mu$ of $\mathcal{A}^{(n)}$ given by $A(\mu_1,\mu_2) = (n-3-\mu_1-\mu_2,\mu_1)$.
This induces the orbifold invariant $Z_{\mathcal{D}^{(n)}}$, which is treated in \cite{evans/pugh:2009v}, and hence its conjugate orbifold invariant $Z_{\mathcal{D}^{(n)\ast}} = Z_{\mathcal{D}^{(n)}} C$, where $C = [\delta_{\lambda,\overline{\lambda}}]$ is the conjugate modular invariant. The conjugate orbifold invariant is given by
\begin{eqnarray*}
Z_{\mathcal{D}^{(3k)\ast}} & = & \frac{1}{3} \sum_{\stackrel{\mu \in P^{(3k)}_+}{\scriptscriptstyle{\mu_1 - \mu_2 \equiv 0 \mathrm{mod}3}}} (\chi_{\mu} + \chi_{A \mu} + \chi_{A^2 \mu}) (\chi_{\overline{\mu}}^{\ast} + \chi_{\overline{A \mu}}^{\ast} + \chi_{\overline{A^2 \mu}}^{\ast}), \qquad k \geq 2, \\
Z_{\mathcal{D}^{(n)\ast}} & = & \sum_{\mu \in P^{(n)}_+} \chi_{\mu} \chi_{\overline{A^{(n-3)(\mu_1-\mu_2)} \mu}}^{\ast}, \qquad n \geq 5, \; n \not \equiv 0 \mathrm{mod} 3.
\end{eqnarray*}
The exponents of $Z_{\mathcal{D}^{(n)\ast}}$ are
$$\mathrm{Exp} = \{ A^k(\lambda_1,\lambda_1) | \; \lambda_1 = 0,1,\ldots, \lfloor (n-3)/2 \rfloor; \; k =0,1,2 \}.$$

It was shown in \cite{evans/pugh:2009ii} that this modular invariant is realised by a braided subfactor with nimrep $\mathcal{D}^{(n)\ast}$ \cite[Figure 12]{evans/pugh:2009i}.
Then as in (\ref{eqn:theta->lambda}), with $\theta_1 = (\lambda_1 + 2\lambda_2 + 3)/24$, $\theta_2 = (2\lambda_1 + \lambda_2 + 3)/24$, we have for each $\lambda_1 = 0,1,\ldots, \lfloor (n-3)/2 \rfloor$:
\begin{center}
\begin{tabular}{|c|c|c|} \hline
$\lambda \in \mathrm{Exp}$ & $(\theta_1,\theta_2) \in [0,1]^2$ & $|\psi^{\lambda}_{\ast}|^2$ \\
\hline $(\lambda_1,\lambda_1)$ & $\left(\frac{\lambda_1+1}{n},\frac{\lambda_1+1}{n}\right)$ & \\
\cline{1-2} $(n-2\lambda_1-3,\lambda_1)$ & $\left(\frac{1}{3},\frac{2}{3}-\frac{\lambda_1+1}{n}\right)$ & $\frac{4}{n} \sin^2\left(\frac{2\pi(\lambda_1+1)}{n}\right)$ \\
\cline{1-2} $(\lambda_1,n-2\lambda_1-3)$ & $\left(\frac{2}{3}-\frac{\lambda_1+1}{n},\frac{1}{3}\right)$ & \\
\hline
\end{tabular}
\end{center}

From (\ref{eqn:moments_SU(3)}),
\begin{equation} \label{eqn:moments-D(n)star}
\int_{\mathbb{T}^2} R_{m,n}(\omega_1,\omega_2) \mathrm{d}\varepsilon(\omega_1,\omega_2) = \frac{1}{6} \sum_{g \in S_3} \sum_{\lambda \in \mathrm{Exp}} (\beta^{(g(\lambda))})^m (\overline{\beta^{(g(\lambda))}})^n |\psi^{g(\lambda)}_{\ast}|^2,
\end{equation}
where $g(\lambda)$ is uniquely given by the pair $(\lambda_1',\lambda_2')$ such that if $\theta_1' = (\lambda_1' + 2\lambda_2' + 3)/3n$ and $\theta_2' = (2\lambda_1' + \lambda_2' + 3)/3n$, then $(\theta_1',\theta_2') = g(\theta_1,\theta_2)$.
For each $j = \lambda_1 + 1 = 1,2,\ldots, \lfloor (n-1)/2 \rfloor$, the $S_3$-orbit of $(e^{2 \pi i \theta_1},e^{2 \pi i \theta_2}) \in \mathbb{T}^2$ under $S_3$ for $(\theta_1,\theta_2) = (j/n,j/n), (1/3,2/3-j/n), (2/3-j/n,1/3)$ give the measure $\mathrm{d}^{((n/j))}$. Then we obtain the following result:

\begin{Thm} \label{Thm:spec_measure-D(n)star}
The spectral measure of $\mathcal{D}^{(n)\ast}$, $n \geq 5$, (over $\mathbb{T}^2$) is
\begin{equation}
\mathrm{d}\varepsilon(\omega_1,\omega_2) = \frac{12}{n} \sum_{j=1}^{\lfloor (n-1)/2 \rfloor} \sin^2(2\pi j/n) \; \mathrm{d}^{((n/j))}(\omega_1,\omega_2),
\end{equation}
where $\mathrm{d}^{((n/j))}$ is defined in Definition \ref{def:4measures}.
\end{Thm}

\subsection{Graph $\mathcal{E}^{(8)}$: $SU(3)_5 \subset SU(6)_1$}

The modular invariant realised by the inclusion $SU(3)_5 \subset SU(6)_1$ \cite{xu:1998, bockenhauer/evans:1999i, evans/pugh:2009ii} is
\begin{eqnarray*}
Z_{\mathcal{E}^{(8)}} & = & |\chi_{(0,0)}+\chi_{(2,2)}|^2 + |\chi_{(0,2)}+\chi_{(3,2)}|^2 + |\chi_{(2,0)}+ \chi_{(2,3)}|^2 + |\chi_{(2,1)}+\chi_{(0,5)}|^2 \\
& & \;\; + |\chi_{(3,0)}+\chi_{(0,3)}|^2 + |\chi_{(1,2)}+\chi_{(5,0)}|^2
\end{eqnarray*}
with exponents
$$\mathrm{Exp} = \{ (0,0), (2,0), (0,2), (3,0), (2,1), (1,2), (0,3), (2,2), (5,0), (3,2), (2,3), (0,5) \}.$$

The inclusion $SU(3)_9 \subset (\mathrm{E}_6)_1$ produces the nimrep $\mathcal{E}^{(8)}$ \cite{xu:1998, bockenhauer/evans:1999i, evans/pugh:2009ii}.
Then as in (\ref{eqn:theta->lambda}), with $\theta_1 = (\lambda_1 + 2\lambda_2 + 3)/24$, $\theta_2 = (2\lambda_1 + \lambda_2 + 3)/24$, for $\lambda = (\lambda_1,\lambda_2) \in \mathrm{Exp}$, we have:
\begin{center}
\begin{tabular}{|c|c|c|} \hline
$\lambda \in \mathrm{Exp}$ & $(\theta_1,\theta_2) \in [0,1]^2$ & $|\psi^{\lambda}_{\ast}|^2$ \\
\hline $(0,0)$, $(5,0)$, $(0,5)$ & $\left(\frac{1}{8},\frac{1}{8}\right)$, $\left(\frac{1}{3},\frac{13}{24}\right)$, $\left(\frac{13}{24},\frac{1}{3}\right)$ & $\frac{2-\sqrt{2}}{24}$ \\
\hline $(2,2)$, $(2,1)$, $(1,2)$ & $\left(\frac{3}{8},\frac{3}{8}\right)$, $\left(\frac{7}{24},\frac{1}{3}\right)$, $\left(\frac{1}{3},\frac{7}{24}\right)$ & $\frac{2+\sqrt{2}}{24}$ \\
\hline $(3,0)$, $(2,3)$, $(0,2)$ & $\left(\frac{1}{4},\frac{3}{8}\right)$, $\left(\frac{11}{24},\frac{5}{12}\right)$, $\left(\frac{7}{24},\frac{5}{24}\right)$ & $\frac{1}{12}$ \\
\hline $(0,3)$, $(3,2)$, $(2,0)$ & $\left(\frac{3}{8},\frac{1}{4}\right)$, $\left(\frac{5}{12},\frac{11}{24}\right)$, $\left(\frac{5}{24},\frac{7}{24}\right)$ & $\frac{1}{12}$ \\
\hline
\end{tabular}
\end{center}

Again, from (\ref{eqn:moments_SU(3)}),
\begin{equation} \label{eqn:moments-E(8)}
\int_{\mathbb{T}^2} R_{m,n}(\omega_1,\omega_2) \mathrm{d}\varepsilon(\omega_1,\omega_2) = \frac{1}{6} \sum_{g \in S_3} \sum_{\lambda \in \mathrm{Exp}} (\beta^{(g(\lambda))})^m (\overline{\beta^{(g(\lambda))}})^n |\psi^{g(\lambda)}_{\ast}|^2.
\end{equation}

\begin{figure}[bt]
\begin{center}
  \includegraphics[width=55mm]{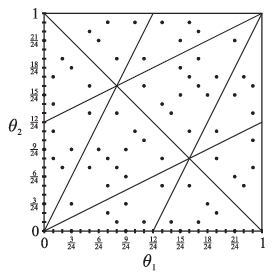}\\
 \caption{The points $(\theta_1,\theta_2)$ given by $g(\lambda)$, where $g \in S_3$, $\lambda \in \mathrm{Exp}$, for $\mathcal{E}^{(8)}$.} \label{fig:poly-11}
\end{center}
\end{figure}

The pairs $(\theta_1,\theta_2)$ given by $g(\lambda)$ for $\lambda \in \mathrm{Exp}$, $g \in S_3$, are illustrated in Figure \ref{fig:poly-11}.
The measure $\mathrm{d}^{((8))}$ is the uniform measure on the $S_3$-orbit of $(e^{2 \pi i \theta_1},e^{2 \pi i \theta_2}) \in \mathbb{T}^2$ under $S_3$ for $(\theta_1,\theta_2) = (1/8,1/8), (1/3,13/24), (13/24,1/3)$, and similarly $\mathrm{d}^{((8/3))}$ is the measure for the $S_3$-orbit for $(\theta_1,\theta_2) = (3/8,3/8), (1/3,7/24), (7/24,1/3)$. The remaining points $(\theta_1,\theta_2)$ that appear in (\ref{eqn:moments-E(8)}) give the measure $\mathrm{d}^{(24/5,1/12)}$. Then the spectral measure for $\mathcal{E}^{(8)}$ is
$$\mathrm{d}\varepsilon = \frac{1}{6} \left( 18 \; \frac{2-\sqrt{2}}{24} \; \mathrm{d}^{((8))} + 18 \; \frac{2+\sqrt{2}}{24} \; \mathrm{d}^{((8/3))} + \frac{36}{12} \; \mathrm{d}^{(24/5,1/12)} \right),$$
and we have obtained the following result:

\begin{Thm} \label{Thm:spec_measure-E(8)}
The spectral measure of $\mathcal{E}^{(8)}$ (over $\mathbb{T}^2$) is
\begin{equation}
\mathrm{d}\varepsilon = \frac{2-\sqrt{2}}{8} \; \mathrm{d}^{((8))} + \frac{2+\sqrt{2}}{8} \; \mathrm{d}^{((8/3))} + \frac{1}{2} \; \mathrm{d}^{(24/5,1/12)},
\end{equation}
where $\mathrm{d}^{((n))}$, $\mathrm{d}^{(n,k)}$ are as in Definition \ref{def:4measures}.
\end{Thm}

\subsection{Graph $\mathcal{E}_1^{(12)}$: $SU(3)_9 \subset (\mathrm{E}_6)_1$}

The modular invariant realised by the inclusion $SU(3)_9 \subset (\mathrm{E}_6)_1$ \cite{xu:1998, bockenhauer/evans:1999ii, evans/pugh:2009ii} is
\begin{eqnarray}
Z_{\mathcal{E}^{(12)}} \;\; = \;\; Z_{\mathcal{E}^{(12)}} C & = & |\chi_{(0,0)}+\chi_{(0,9)}+\chi_{(9,0)}+\chi_{(4,4)}+\chi_{(4,1)}+\chi_{(1,4)}|^2 \nonumber \\
& & \;\; + 2 |\chi_{(2,2)}+\chi_{(2,5)}+\chi_{(5,2)}|^2, \label{eqn:Z(E(12))}
\end{eqnarray}
and its exponents are
$$\mathrm{Exp} = \{ (0,0), (4,1), (1,4), (4,4), (9,0), (0,9), \textrm{ and twice } (2,2), (5,2), (2,5) \}.$$

The inclusion $SU(3)_9 \subset (\mathrm{E}_6)_1$ produces the nimrep $\mathcal{E}_1^{(12)}$ \cite{xu:1998, bockenhauer/evans:1999ii, evans/pugh:2009ii}.
Then as in (\ref{eqn:theta->lambda}), with $\theta_1 = (\lambda_1 + 2\lambda_2 + 3)/24$, $\theta_2 = (2\lambda_1 + \lambda_2 + 3)/24$, for $\lambda = (\lambda_1,\lambda_2) \in \mathrm{Exp}$, we have:
\begin{center}
\begin{tabular}{|c|c|c|} \hline
$\lambda \in \mathrm{Exp}$ & $(\theta_1,\theta_2) \in [0,1]^2$ & $|\psi^{\lambda}_{\ast}|^2$ \\
\hline $(0,0)$, $(9,0)$, $(0,9)$ & $\left(\frac{1}{12},\frac{1}{12}\right)$, $\left(\frac{7}{12},\frac{1}{3}\right)$, $\left(\frac{1}{3},\frac{7}{12}\right)$ & $\frac{2-\sqrt{3}}{36}$ \\
\hline $(4,4)$, $(4,1)$, $(1,4)$ & $\left(\frac{5}{12},\frac{5}{12}\right)$, $\left(\frac{1}{3},\frac{1}{4}\right)$, $\left(\frac{1}{4},\frac{1}{3}\right)$ & $\frac{2+\sqrt{3}}{36}$ \\
\hline $(2,2)$, $(5,2)$, $(2,5)$ & $\left(\frac{1}{4},\frac{1}{4}\right)$, $\left(\frac{5}{12},\frac{1}{3}\right)$, $\left(\frac{1}{3},\frac{5}{12}\right)$ & $\frac{2}{9}$ \\
\hline
\end{tabular}
\end{center}

\begin{figure}[bt]
\begin{center}
  \includegraphics[width=55mm]{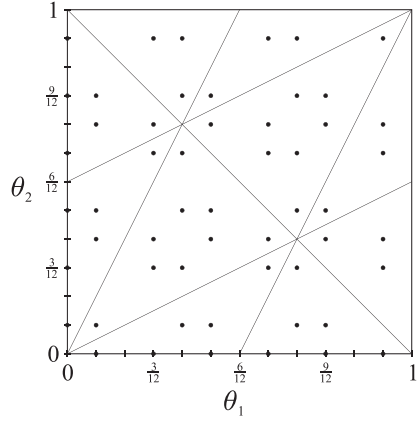}\\
 \caption{The points $(\theta_1,\theta_2) \in \{ g(\lambda) | \; \lambda \in \mathrm{Exp}, g \in S_3 \}$ for $\mathcal{E}_1^{(12)}$.} \label{fig:poly-12}
\end{center}
\end{figure}

We illustrate the pairs $(\theta_1,\theta_2)$ given by $g(\lambda)$ for $\lambda \in \mathrm{Exp}$, $g \in S_3$, in Figure \ref{fig:poly-12}.
We obtain the following spectral measure for $\mathcal{E}_1^{(12)}$:
\begin{eqnarray}
\mathrm{d}\varepsilon & = & \frac{1}{6} \left( 18 \; \frac{2-\sqrt{3}}{36} \; \mathrm{d}^{((12))} + 18 \; \frac{2+\sqrt{3}}{36} \; \mathrm{d}^{((12/5))} + 18 \; \frac{2}{9} \; \mathrm{d}^{((4))} \right) \nonumber \\
& = & \frac{2-\sqrt{3}}{12} \; \mathrm{d}^{((12))} + \frac{2+\sqrt{3}}{12} \; \mathrm{d}^{((12/5))} + \frac{2}{3} \; \mathrm{d}^{((4))}.
\end{eqnarray}
Now for the pairs $(\theta_1,\theta_2)$ given by the measure $\mathrm{d}^{((4))}$, we have $J(\theta_1,\theta_2)^2 = 64\pi^4$, so that the spectral measure for $\mathcal{E}_1^{(12)}$ is
$$\mathrm{d}\varepsilon = \frac{2-\sqrt{3}}{12} \; \mathrm{d}^{((12))} + \frac{2+\sqrt{3}}{12} \; \mathrm{d}^{((12/5))} + \frac{1}{96\pi^4} J^2 \; \mathrm{d}^{((4))},$$
and using (\ref{eqn:measure-relation2}) we obtain the following result:

\begin{Thm} \label{Thm:spec_measure-E1(12)}
The spectral measure of $\mathcal{E}_1^{(12)}$ (over $\mathbb{T}^2$) is
\begin{equation}
\mathrm{d}\varepsilon = \frac{2-\sqrt{3}}{12} \; \mathrm{d}^{((12))} + \frac{2+\sqrt{3}}{12} \; \mathrm{d}^{((12/5))} + \frac{1}{36\pi^4} J^2 \; \mathrm{d}^{(4)},
\end{equation}
where $\mathrm{d}^{(n)}$, $\mathrm{d}^{((n))}$ are as in Definition \ref{def:4measures}.
\end{Thm}

\subsection{Graph $\mathcal{E}_2^{(12)}$: $SU(3)_9 \subset (\mathrm{E}_6)_1 \rtimes \mathbb{Z}_3$}

The modular invariant realised by the inclusion $SU(3)_9 \subset (\mathrm{E}_6)_1 \rtimes \mathbb{Z}_3$ is again the modular invariant $Z_{\mathcal{E}^{(12)}}$ given in (\ref{eqn:Z(E(12))}). However the nimrep obtained from the inclusion $SU(3)_9 \subset (\mathrm{E}_6)_1 \rtimes \mathbb{Z}_3$ is the graph $\mathcal{E}_2^{(12)}$ \cite{bockenhauer/evans:2002, evans/pugh:2009ii}. Hence the graphs $\mathcal{E}_2^{(12)}$ and $\mathcal{E}_1^{(12)}$ are isospectral.
Then as in (\ref{eqn:theta->lambda}), with $\theta_1 = (\lambda_1 + 2\lambda_2 + 3)/24$, $\theta_2 = (2\lambda_1 + \lambda_2 + 3)/24$, for $\lambda = (\lambda_1,\lambda_2) \in \mathrm{Exp}$, we have:
\begin{center}
\begin{tabular}{|c|c|c|} \hline
$\lambda \in \mathrm{Exp}$ & $(\theta_1,\theta_2) \in [0,1]^2$ & $|\psi^{\lambda}_{\ast}|^2$ \\
\hline $(0,0)$, $(9,0)$, $(0,9)$ & $\left(\frac{1}{12},\frac{1}{12}\right)$, $\left(\frac{7}{12},\frac{1}{3}\right)$, $\left(\frac{1}{3},\frac{7}{12}\right)$ & $\frac{2+\sqrt{3}}{36}$ \\
\hline $(4,4)$, $(4,1)$, $(1,4)$ & $\left(\frac{5}{12},\frac{5}{12}\right)$, $\left(\frac{1}{3},\frac{1}{4}\right)$, $\left(\frac{1}{4},\frac{1}{3}\right)$ & $\frac{2-\sqrt{3}}{36}$ \\
\hline $(2,2)$, $(5,2)$, $(2,5)$ & $\left(\frac{1}{4},\frac{1}{4}\right)$, $\left(\frac{5}{12},\frac{1}{3}\right)$, $\left(\frac{1}{3},\frac{5}{12}\right)$ & $\frac{2}{9}$ \\
\hline
\end{tabular}
\end{center}

We see that the spectral measure for $\mathcal{E}_2^{(12)}$ is identical to that for $\mathcal{E}_1^{(12)}$, except that the weights $(2+\sqrt{3})/12$ and $(2-\sqrt{3})/12$ are interchanged, giving:

\begin{Thm} \label{Thm:spec_measure-E2(12)}
The spectral measure of $\mathcal{E}_2^{(12)}$ (over $\mathbb{T}^2$) is
\begin{equation}
\mathrm{d}\varepsilon = \frac{2+\sqrt{3}}{12} \; \mathrm{d}^{((12))} + \frac{2-\sqrt{3}}{12} \; \mathrm{d}^{((12/5))} + \frac{1}{36\pi^4} J^2 \; \mathrm{d}^{(4)},
\end{equation}
where $\mathrm{d}^{(n)}$, $\mathrm{d}^{((n))}$ are as in Definition \ref{def:4measures}.
\end{Thm}

\subsection{Graph $\mathcal{E}_5^{(12)}$}

The Moore-Seiberg invariant
\begin{eqnarray}
Z_{\mathcal{E}_{MS}^{(12)}} & = & |\chi_{(0,0)}+\chi_{(0,9)}+\chi_{(9,0)}|^2 + |\chi_{(2,2)}+\chi_{(2,5)}+\chi_{(5,2)}|^2 + 2 |\chi_{(3,3)}|^2 \nonumber \\
& & \;\; + |\chi_{(0,3)}+\chi_{(6,0)}+ \chi_{(3,6)}|^2 + |\chi_{(3,0)}+\chi_{(0,6)}+\chi_{(6,3)}|^2 + |\chi_{(4,4)}+\chi_{(4,1)}+\chi_{(1,4)}|^2 \nonumber \\
& & \;\; + (\chi_{(1,1)}+\chi_{(1,7)}+\chi_{(7,1)})\chi_{(3,3)}^{\ast} + \chi_{(3,3)}(\chi_{(1,1)}^{\ast}+\chi_{(1,7)}^{\ast}+\chi_{(7,1)}^{\ast}).
\end{eqnarray}
has exponents
\begin{eqnarray*}
\mathrm{Exp} & = & \{ (0,0), (3,0), (0,3), (2,2), (4,1), (1,4), (6,0), (0,6), (5,2), (2,5), \\
& & \quad (4,4), (9,0), (6,3), (3,6), (0,9) \textrm{ and twice } (3,3) \}.
\end{eqnarray*}
It is realised by a braided subfactor which produces the nimrep $\mathcal{E}_5^{(12)}$ \cite[Section 5.4]{evans/pugh:2009ii}.
Then as in (\ref{eqn:theta->lambda}), with $\theta_1 = (\lambda_1 + 2\lambda_2 + 3)/24$, $\theta_2 = (2\lambda_1 + \lambda_2 + 3)/24$, for $\lambda = (\lambda_1,\lambda_2) \in \mathrm{Exp}$, we have:
\begin{center}
\begin{tabular}{|c|c|c|} \hline
$\lambda \in \mathrm{Exp}$ & $(\theta_1,\theta_2) \in [0,1]^2$ & $|\psi^{\lambda}_{\ast}|^2$ \\
\hline $(0,0)$, $(9,0)$, $(0,9)$ & $\left(\frac{1}{12},\frac{1}{12}\right)$, $\left(\frac{7}{12},\frac{1}{3}\right)$, $\left(\frac{1}{3},\frac{7}{12}\right)$ & $\frac{1}{36}$ \\
\hline $(4,4)$, $(4,1)$, $(1,4)$ & $\left(\frac{5}{12},\frac{5}{12}\right)$, $\left(\frac{1}{3},\frac{1}{4}\right)$, $\left(\frac{1}{4},\frac{1}{3}\right)$ & $\frac{1}{36}$ \\
\hline $(2,2)$, $(5,2)$, $(2,5)$ & $\left(\frac{1}{4},\frac{1}{4}\right)$, $\left(\frac{5}{12},\frac{1}{3}\right)$, $\left(\frac{1}{3},\frac{5}{12}\right)$ & $\frac{1}{9}$ \\
\hline $(3,0)$, $(6,3)$, $(0,6)$ & $\left(\frac{1}{6},\frac{1}{4}\right)$, $\left(\frac{5}{12},\frac{1}{2}\right)$, $\left(\frac{5}{12},\frac{1}{4}\right)$ & 0 \\
$(0,3)$, $(3,6)$, $(6,0)$ & $\left(\frac{1}{4},\frac{1}{6}\right)$, $\left(\frac{1}{2},\frac{5}{12}\right)$, $\left(\frac{1}{4},\frac{5}{12}\right)$ & \\
\hline $(3,3)$ & $\left(\frac{1}{3},\frac{1}{3}\right)$ & $\frac{1}{2}$ \\
\hline
\end{tabular}
\end{center}
Then we obtain the following spectral measure for $\mathcal{E}_5^{(12)}$:
\begin{eqnarray}
\mathrm{d}\varepsilon & = & \frac{1}{6} \left( \frac{18}{9} \; \mathrm{d}^{((12))} + \frac{18}{9} \; \mathrm{d}^{((12/5))} + \frac{18}{9} \; \mathrm{d}^{((4))} \right) \nonumber \\
& = & \frac{1}{3} \left( \mathrm{d}^{((12))} + \mathrm{d}^{((12/5))} + \mathrm{d}^{((4))} \right), \nonumber
\end{eqnarray}
and using (\ref{eqn:measure-relation2}) we obtain the following result:

\begin{Thm} \label{Thm:spec_measure-E5(12)}
The spectral measure of $\mathcal{E}_5^{(12)}$ (over $\mathbb{T}^2$) is
\begin{equation}
\mathrm{d}\varepsilon = \frac{1}{3} \left( \mathrm{d}^{((12))} + \mathrm{d}^{((12/5))} + \frac{1}{24\pi^4} J^2 \; \mathrm{d}^{(4)} \right),
\end{equation}
where $\mathrm{d}^{(n)}$, $\mathrm{d}^{((n))}$ are as in Definition \ref{def:4measures}.
\end{Thm}

\subsection{Graph $\mathcal{E}_4^{(12)}$} \label{sec:E4(12)}

It has not yet been shown that the graph $\mathcal{E}_4^{(12)}$ is a nimrep obtained from an inclusion which realises the conjugate Moore-Seiberg modular invariant $Z_{\mathcal{E}_{MS}^{(12)\ast}} = Z_{\mathcal{E}_{MS}^{(12)}} C$, given by
\begin{eqnarray}
Z_{\mathcal{E}_{MS}^{(12)\ast}} & = & |\chi_{(0,0)}+\chi_{(0,9)}+\chi_{(9,0)}|^2 + |\chi_{(2,2)}+\chi_{(2,5)}+\chi_{(5,2)}|^2 + 2 |\chi_{(3,3)}|^2 \nonumber \\
& & \;\; + (\chi_{(0,3)}+\chi_{(6,0)}+ \chi_{(3,6)})(\chi_{(3,0)}^{\ast}+\chi_{(0,6)}^{\ast}+\chi_{(6,3)}^{\ast}) \nonumber \\
& & \;\; + (\chi_{(3,0)}+\chi_{(0,6)}+\chi_{(6,3)})(\chi_{(0,3)}^{\ast}+\chi_{(6,0)}^{\ast}+ \chi_{(3,6)}^{\ast}) + |\chi_{(4,4)}+\chi_{(4,1)}+\chi_{(1,4)}|^2 \nonumber \\
& & \;\; + (\chi_{(1,1)}+\chi_{(1,7)}+\chi_{(7,1)})\chi_{(3,3)}^{\ast} + \chi_{(3,3)}(\chi_{(1,1)}^{\ast}+\chi_{(1,7)}^{\ast}+\chi_{(7,1)}^{\ast}).
\end{eqnarray}
However, it is known that $\mathcal{E}_4^{(12)}$ is a nimrep \cite{di_francesco/zuber:1990}, and it can be checked by hand that the eigenvalues of $\mathcal{E}_4^{(12)}$ are given by $S_{\mu, \rho}/S_{\mu,1}$, where $\mu$ runs over the exponents of $Z_{\mathcal{E}_{MS}^{(12)\ast}}$:
$$\mathrm{Exp} = \{ (0,0), (2,2), (4,1), (1,4), (5,2), (2,5), (4,4), (9,0), (0,9), \textrm{ and twice } (3,3) \}.$$
With $\theta_1 = (\lambda_1 + 2\lambda_2 + 3)/24$, $\theta_2 = (2\lambda_1 + \lambda_2 + 3)/24$, for $\lambda = (\lambda_1,\lambda_2) \in \mathrm{Exp}$, we have:
\begin{center}
\begin{tabular}{|c|c|c|} \hline
$\lambda \in \mathrm{Exp}$ & $(\theta_1,\theta_2) \in [0,1]^2$ & $|\psi^{\lambda}_{\ast}|^2$ \\
\hline $(0,0)$, $(9,0)$, $(0,9)$ & $\left(\frac{1}{12},\frac{1}{12}\right)$, $\left(\frac{7}{12},\frac{1}{3}\right)$, $\left(\frac{1}{3},\frac{7}{12}\right)$ & $\frac{1}{36}$ \\
\hline $(4,4)$, $(4,1)$, $(1,4)$ & $\left(\frac{5}{12},\frac{5}{12}\right)$, $\left(\frac{1}{3},\frac{1}{4}\right)$, $\left(\frac{1}{4},\frac{1}{3}\right)$ & $\frac{1}{36}$ \\
\hline $(2,2)$, $(5,2)$, $(2,5)$ & $\left(\frac{1}{4},\frac{1}{4}\right)$, $\left(\frac{5}{12},\frac{1}{3}\right)$, $\left(\frac{1}{3},\frac{5}{12}\right)$ & $\frac{1}{9}$ \\
\hline $(3,3)$ & $\left(\frac{1}{3},\frac{1}{3}\right)$ & $\frac{1}{2}$ \\
\hline
\end{tabular}
\end{center}
The pairs $(\theta_1,\theta_2)$ given by $g(\lambda)$ for $\lambda \in \mathrm{Exp} \setminus \{(3,3)\}$, $g \in S_3$, have all appeared in the computations of the spectral measures for the graphs $\mathcal{E}_1^{(12)}$, $\mathcal{E}_2^{(12)}$, hence the measures which give these points are known. For the remaining points $g((3,3))$, $g \in S_3$, we have $J (g(e^{2 \pi i/3}, e^{2 \pi i/3}))^2 = 108 \pi^4$, whilst $J (e^{2 \pi i \theta_1}, e^{2 \pi i \theta_2})^2 = 0$ for the other points that are given by the measure $\mathrm{d}^{(3)}$, since they map to the boundary of the discoid $\mathfrak{D}$.
Then we obtain the following spectral measure for $\mathcal{E}_4^{(12)}$:
\begin{eqnarray}
\mathrm{d}\varepsilon & = & \frac{1}{6} \left( \frac{18}{36} \; \mathrm{d}^{((12))} + \frac{18}{36} \; \mathrm{d}^{((12/5))} + \frac{18}{9} \; \mathrm{d}^{((4))} + \frac{27}{2} \; \frac{1}{108\pi^4} J^2 \; \mathrm{d}^{(3)} \right) \nonumber \\
& = & \frac{1}{12} \; \mathrm{d}^{((12))} + \frac{1}{12} \; \mathrm{d}^{((12/5))} + 2 \mathrm{d}^{((4))} + \frac{1}{8\pi^4} J^2 \; \mathrm{d}^{(3)}, \nonumber
\end{eqnarray}
and using (\ref{eqn:measure-relation2}) we obtain the following result:

\begin{Thm} \label{Thm:spec_measure-E4(12)}
The spectral measure of $\mathcal{E}_4^{(12)}$ (over $\mathbb{T}^2$) is
\begin{equation}
\mathrm{d}\varepsilon = \frac{1}{12} \; \mathrm{d}^{((12))} + \frac{1}{12} \; \mathrm{d}^{((12/5))} + \frac{1}{12\pi^4} J^2 \; \mathrm{d}^{(4)} + \frac{1}{8\pi^4} J^2 \; \mathrm{d}^{(3)},
\end{equation}
where $\mathrm{d}^{(n)}$, $\mathrm{d}^{((n))}$ are as in Definition \ref{def:4measures}.
\end{Thm}

\subsection{Graph $\mathcal{E}^{(24)}$: $SU(3)_{21} \subset (E_7)_1$}

The modular invariant realised by the inclusion $SU(3)_{21} \subset (E_7)_1$ \cite{xu:1998, bockenhauer/evans:1999ii, evans/pugh:2009ii} is
\begin{eqnarray}
Z_{\mathcal{E}^{(24)}} & = & |\chi_{(0,0)}+\chi_{(4,4)}+\chi_{(6,6)}+\chi_{(10,10)}+\chi_{(21,0)}+\chi_{(0,21)}+\chi_{(13,4)}+\chi_{(4,13)} \nonumber \\
& & \;\; +\chi_{(10,1)}+\chi_{(1,10)}+\chi_{(9,6)}+\chi_{(6,9)}|^2 \nonumber \\
& & \;\; + |\chi_{(15,6)}+\chi_{(6,15)}+\chi_{(15,0)}+\chi_{(0,15)}+\chi_{(10,7)}+\chi_{(7,10)}+\chi_{(10,4)} \nonumber \\
& & \;\; \quad +\chi_{(4,10)}+\chi_{(7,4)}+\chi_{(4,7)}+\chi_{(6,0)}+\chi_{(0,6)}|^2.
\end{eqnarray}
with exponents
\begin{eqnarray*}
\mathrm{Exp} & = & \{ (0,0), (6,0), (0,6), (10,1), (7,4), (4,7), (1,10), (6,6), (10,4), (4,10), \\
& & \quad (15,0), (9,6), (6,9), (0,15), (13,4), (10,7), (7,10), (4,13), (10,10), \\
& & \quad (21,0), (15,6), (6,15), (0,21) \}.
\end{eqnarray*}

The inclusion $SU(3)_{21} \subset (E_7)_1$ produces the nimrep $\mathcal{E}^{(24)}$ \cite{xu:1998, bockenhauer/evans:1999ii, evans/pugh:2009ii}.
Then as in (\ref{eqn:theta->lambda}), with $\theta_1 = (\lambda_1 + 2\lambda_2 + 3)/24$, $\theta_2 = (2\lambda_1 + \lambda_2 + 3)/24$, for $\lambda = (\lambda_1,\lambda_2) \in \mathrm{Exp}$, we have:
\begin{center}
\begin{tabular}{|c|c|c|} \hline
$\lambda \in \mathrm{Exp}$ & $(\theta_1,\theta_2) \in [0,1]^2$ & $|\psi^{\lambda}_{\ast}|^2$ \\
\hline $(0,0)$, $(21,0)$, $(0,21)$ & $\left(\frac{1}{24},\frac{1}{24}\right)$, $\left(\frac{5}{8},\frac{1}{3}\right)$, $\left(\frac{1}{3},\frac{5}{8}\right)$ & $\frac{6-2\sqrt{3}-\sqrt{6}}{144}$ \\
\hline $(4,4)$, $(13,4)$, $(4,13)$ & $\left(\frac{5}{24},\frac{5}{24}\right)$, $\left(\frac{11}{24},\frac{1}{3}\right)$, $\left(\frac{1}{3},\frac{11}{24}\right)$ & $\frac{6+2\sqrt{3}-\sqrt{6}}{144}$ \\
\hline $(6,6)$, $(9,6)$, $(6,9)$ & $\left(\frac{7}{24},\frac{7}{24}\right)$, $\left(\frac{3}{8},\frac{1}{3}\right)$, $\left(\frac{1}{3},\frac{3}{8}\right)$ & $\frac{6+2\sqrt{3}+\sqrt{6}}{144}$ \\
\hline $(10,10)$, $(10,1)$, $(1,10)$ & $\left(\frac{11}{24},\frac{11}{24}\right)$, $\left(\frac{5}{24},\frac{1}{3}\right)$, $\left(\frac{1}{3},\frac{5}{24}\right)$ & $\frac{6-2\sqrt{3}+\sqrt{6}}{144}$ \\
\hline $(6,0)$, $(15,6)$, $(0,15)$ & $\left(\frac{5}{24},\frac{1}{8}\right)$, $\left(\frac{13}{24},\frac{5}{12}\right)$, $\left(\frac{1}{4},\frac{11}{24}\right)$ & $\frac{2-\sqrt{2}}{48}$ \\
$(0,6)$, $(6,15)$, $(15,0)$ & $\left(\frac{1}{8},\frac{5}{24}\right)$, $\left(\frac{5}{12},\frac{13}{24}\right)$, $\left(\frac{11}{24},\frac{1}{4}\right)$ & \\
\hline $(7,4)$, $(10,7)$, $(4,10)$ & $\left(\frac{7}{24},\frac{1}{4}\right)$, $\left(\frac{5}{12},\frac{3}{8}\right)$, $\left(\frac{7}{24},\frac{3}{8}\right)$ & $\frac{2+\sqrt{2}}{48}$ \\
$(4,7)$, $(7,10)$, $(10,4)$ & $\left(\frac{1}{4},\frac{7}{24}\right)$, $\left(\frac{3}{8},\frac{5}{12}\right)$, $\left(\frac{3}{8},\frac{7}{24}\right)$ & \\
\hline
\end{tabular}
\end{center}
Then we obtain:

\begin{Thm} \label{Thm:spec_measure-E(24)}
The spectral measure of $\mathcal{E}^{(24)}$ (over $\mathbb{T}^2$) is
\begin{eqnarray}
\mathrm{d}\varepsilon & = & \frac{6-2\sqrt{3}-\sqrt{6}}{48} \; \mathrm{d}^{((24))} + \frac{6+2\sqrt{3}-\sqrt{6}}{48} \; \mathrm{d}^{((24/5))} + \frac{6+2\sqrt{3}+\sqrt{6}}{48} \; \mathrm{d}^{((24/7))} \nonumber \\
& & + \frac{6-2\sqrt{3}+\sqrt{6}}{48} \; \mathrm{d}^{((24/11))} + \frac{2-\sqrt{2}}{8} \; \mathrm{d}^{(8,1/12)} + \frac{2+\sqrt{2}}{8} \; \mathrm{d}^{(4,1/24)},
\end{eqnarray}
where $\mathrm{d}^{((n))}$, $\mathrm{d}^{(n,k)}$ are as in Definition \ref{def:4measures}.
\end{Thm}

\newpage

\section{Spectral measures for finite subgroups of $SU(3)$} \label{sect:spec_measure-subgroupsSU(3)}

\renewcommand{\arraystretch}{1}

\begin{center}
\begin{table}[tb]
\begin{tabular}{|c|c|c|c|} \hline
$\mathcal{ADE}$ graph & Type & Subgroup $\Gamma \subset SU(3)$ & $|\Gamma|$ \\
\hline\hline ($ADE$) & - & B: finite subgroups of $SU(2) \subset SU(3)$ & - \\
\hline $\mathcal{A}^{(n)}$ & I & A: $\mathbb{Z}_{n-2} \times \mathbb{Z}_{n-2}$ & $(n-2)^2$ \\
\hline - & - & A: $\mathbb{Z}_m \times \mathbb{Z}_n \;\;$ ($m \neq n \neq 3$) & $mn$ \\
\hline $\mathcal{D}^{(n)} \;$ ($n \equiv 0 \textrm{ mod } 3$) & I & C: $\Delta (3(n-3)^2) = (\mathbb{Z}_{n-3} \times \mathbb{Z}_{n-3}) \rtimes \mathbb{Z}_3$ & $3(n-3)^2$ \\
\hline $\mathcal{D}^{(n)} \;$ ($n \not \equiv 0 \textrm{ mod } 3$) & II & - & - \\
\hline - & - & C: $\Delta (3n^2) = (\mathbb{Z}_n \times \mathbb{Z}_n) \rtimes \mathbb{Z}_3$, & $3n^2$ \\
& & ($n \not \equiv 0 \textrm{ mod } 3$) & \\
\hline - & - & D: $\Delta (6n^2) = (\mathbb{Z}_n \times \mathbb{Z}_n) \rtimes S_3$ & $6n^2$ \\
\hline $\mathcal{A}^{(n)\ast}$ & II & - & - \\
\hline $\mathcal{D}^{(n) \ast} \;$ ($n \geq 7$) & II & A: $\mathbb{Z}_{\lfloor (n+1)/2 \rfloor} \times \mathbb{Z}_{3}$ & $3 \lfloor (n+1)/2 \rfloor$ \\
\hline $\mathcal{E}^{(8)}$ & I & E $= \Sigma (36 \times 3) = \Delta (3.3^2) \rtimes \mathbb{Z}_4$ & 108 \\
\hline $\mathcal{E}^{(8)\ast}$ & II & - & - \\
\hline $\mathcal{E}_1^{(12)}$ & I & F $= \Sigma (72 \times 3)$ & 216 \\
\hline $\mathcal{E}_2^{(12)}$ & II & G $= \Sigma (216 \times 3)$ & 648 \\
\hline ($\mathcal{E}_3^{(12)}$) & - & B $\times \mathbb{Z}_3$: $BD_4 \times \mathbb{Z}_3$ & 24 \\
\hline $\mathcal{E}_4^{(12)}$ & (II) & L $= \Sigma (360 \times 3) \cong TA_6$ & 1080 \\
\hline $\mathcal{E}_5^{(12)}$ & II & K $\cong TPSL(2,7)$ & 504 \\
\hline $\mathcal{E}^{(24)}$ & I & - & - \\
\hline - & - & H $= \Sigma (60) \cong A_5$ & 60 \\
\hline - & - & I $= \Sigma (168) \cong PSL(2,7)$ & 168 \\
\hline - & - & J $\cong TA_5$ & 180 \\
\hline
\end{tabular} \\
\caption{Relationship between $\mathcal{ADE}$ graphs and subgroups $\Gamma$ of $SU(3)$.} \label{Table:ADE-subgroupsSU(3)}
\end{table}
\end{center}

The classification of finite subgroups of $SU(3)$ is due to \cite{miller/blichfeldt/dickson:1916,fairbairn/fulton/klink:1964,bovier/luling/wyler:1981, yau/yu:1993}. Clearly any finite subgroup of $SU(2)$ is a finite subgroup of $SU(3)$, since we can embed $SU(2)$ in $SU(3)$ by sending $g \rightarrow g \oplus 1 \in SU(3)$, for any $g \subset SU(2)$. These subgroups of $SU(3)$ are called type B. There are three other infinite series of finite groups, called types A, C, D. The groups of type A are the diagonal abelian groups, which correspond to an embedding of the two torus $\mathbb{T}^2$ in $SU(3)$ given by
$$(\rho|_{\mathbb{T}^2})(\omega_1,\omega_2) = \left( \begin{array}{ccc} \omega_1 & 0 & 0 \\ 0 & \omega_2^{-1} & 0 \\ 0 & 0 & \omega_1^{-1}\omega_2 \end{array} \right),$$
for $(\omega_1,\omega_2) \in \mathbb{T}^2$. The groups of type C are the groups $\Delta (3n^2)$, and those of type D are the groups $\Delta(6n^2)$. These ternary trihedral groups are considered in \cite{bovier/luling/wyler:1981, luhn/nasri/ramond:2007, escobar/luhn:2009} and generalize the binary dihedral subgroups of $SU(2)$. There are also eight exceptional groups E-L. The complete list of finite subgroups of $SU(3)$ is given in Table \ref{Table:ADE-subgroupsSU(3)}. Here Type denotes the type of the inclusion found in \cite{evans/pugh:2009ii} which yielded the $\mathcal{ADE}$ graph as a nimrep. An inclusion $N \subset M$ with dual canonical endomorphism $\theta$ is called type I if and only if one of the following equivalent conditions hold \cite[Proposition 3.2]{bockenhauer/evans:2000}:
\begin{itemize}
\item[1.] $Z_{\lambda,0} = \langle \theta, \lambda \rangle$ for all $\lambda \in {}_N \mathcal{X}_N$.
\item[2.] $Z_{0,\lambda} = \langle \theta, \lambda \rangle$ for all $\lambda \in {}_N \mathcal{X}_N$.
\item[3.] Chiral locality holds: $\varepsilon(\theta,\theta)v^2 = v^2$.
\end{itemize}
Otherwise the inclusion is called type II. (In the context of nets of subfactors $N(I) \subset M(I)$, for $I$ a subinterval of the circle, type I means locality of the extended net $M(I)$). We emphasise that the type of the modular invariant is not well defined, as illustrated by the case of $Z_{\mathcal{E}^{(12)}}$ in Section \ref{sect:intro}.
The $\mathcal{ADE}$ graph $\mathcal{E}_3^{(12)}$ is a nimrep which is isospectral to the graphs $\mathcal{E}_1^{(12)}$ and $\mathcal{E}_2^{(12)}$ \cite{di_francesco/zuber:1990}, however Ocneanu ruled it out as a candidate for the modular invariant $Z_{\mathcal{E}^{(12)}}$ by asserting that it did not support a valid cell system \cite{ocneanu:2002}. This graph was ruled out as a natural candidate in Section 5.2 of \cite{evans:2002}.
Caution should be taken regarding the type for the graph $\mathcal{E}_4^{(12)}$. Although we have not yet shown that $\mathcal{E}_4^{(12)}$ is the nimrep obtained from an inclusion, it was shown in \cite{evans/pugh:2009ii} that such an inclusion would be of type II.

\begin{figure}[tb]
\begin{center}
  \includegraphics[width=70mm]{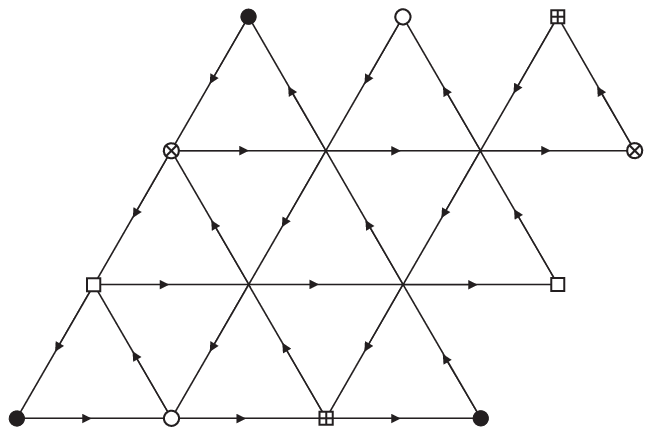}\\
 \caption{$\mathbb{Z}_{n-2} \times \mathbb{Z}_{n-2}$ for $n = 6$; vertices which have the same symbol are identified.}\label{fig:extended-Agraph}
\end{center}
\end{figure}

\begin{figure}[tb]
\begin{minipage}[t]{6cm}
\begin{center}
  \includegraphics[width=50mm]{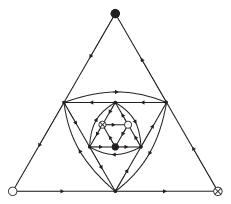}\\
 \caption{$\mathbb{Z}_{p} \times \mathbb{Z}_{3}$ for $p = 3$.}\label{fig:extended-D(star)graph}
\end{center}
\end{minipage}
\hfill
\begin{minipage}[t]{8cm}
\begin{center}
  \includegraphics[width=70mm]{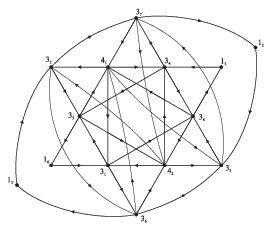}\\
 \caption{E $= \sum (36 \times 3)$}\label{fig:extended-E(8)graph}
\end{center}
\end{minipage}
\end{figure}

\begin{figure}[tb]
\begin{minipage}[t]{6cm}
\begin{center}
  \includegraphics[width=50mm]{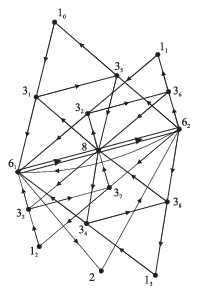}\\
 \caption{F $= \sum (72 \times 3)$}\label{fig:extended-E1(12)graph}
\end{center}
\end{minipage}
\hfill
\begin{minipage}[t]{6cm}
\begin{center}
  \includegraphics[width=50mm]{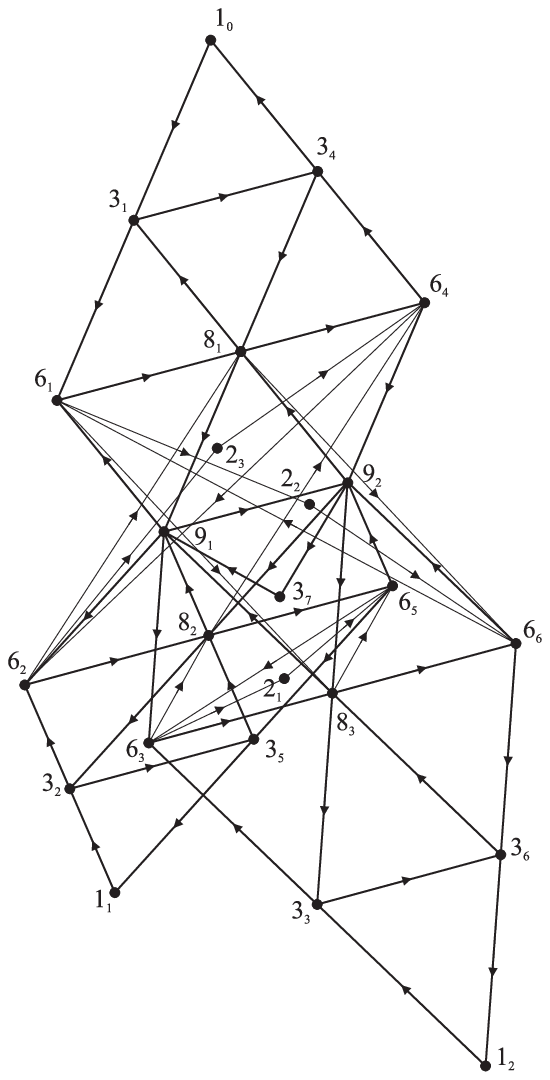}\\
 \caption{G $= \sum (216 \times 3)$}\label{fig:extended-E2(12)graph}
\end{center}
\end{minipage}
\end{figure}

\begin{figure}[tb]
\begin{minipage}[t]{7cm}
\begin{center}
  \includegraphics[width=60mm]{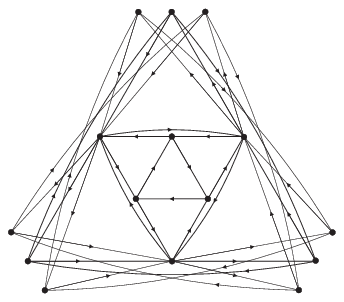}\\
 \caption{$\widehat{\mathcal{D}_4} \otimes \sigma_{123}$}\label{fig:extended-E3(12)graph}
\end{center}
\end{minipage}
\hfill
\begin{minipage}[t]{6cm}
\begin{center}
  \includegraphics[width=50mm]{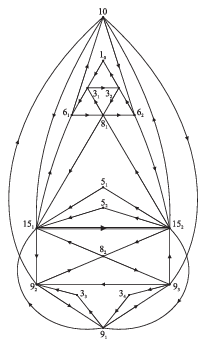}\\
 \caption{L $= \sum (360 \times 3)$}\label{fig:extended-E4(12)graph}
\end{center}
\end{minipage}
\end{figure}

\begin{figure}[tb]
\begin{minipage}[t]{7cm}
\begin{center}
  \includegraphics[width=60mm]{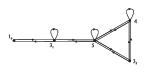}\\
 \caption{H $= \sum (60)$}\label{fig:extended-Hgraph}
\end{center}
\end{minipage}
\hfill
\begin{minipage}[t]{6cm}
\begin{center}
  \includegraphics[width=30mm]{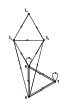}\\
 \caption{I $= \sum (168)$}\label{fig:extended-Igraph}
\end{center}
\end{minipage}
\end{figure}

\begin{figure}[tb]
\begin{minipage}[t]{8cm}
\begin{center}
  \includegraphics[width=75mm]{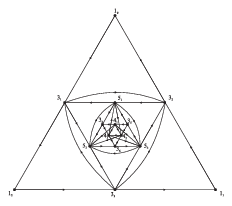}\\
 \caption{J}\label{fig:extended-Jgraph}
\end{center}
\end{minipage}
\hfill
\begin{minipage}[t]{7cm}
\begin{center}
  \includegraphics[width=65mm]{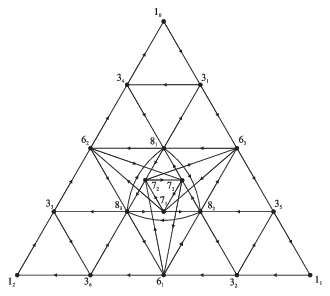}\\
 \caption{K}\label{fig:extended-Kgraph}
\end{center}
\end{minipage}
\end{figure}

The fundamental representation $\rho$ of $SU(3)$ corresponds to the vertex $(1,0)$ of the graph $\mathcal{A}^{(\infty)}$. The McKay graph $\mathcal{G}_{\Gamma}$ is the the fusion graph of $\rho$ acting on the irreducible representations of $\Gamma$.
For most of the graphs $\mathcal{G}_{\Gamma}$ there is a corresponding $SU(3)$ $\mathcal{ADE}$ graph/quiver $\mathcal{G}$, which is obtained from $\mathcal{G}_{\Gamma}$ by now removing more than one vertex, and all the edges that start or end at those vertices, as well as possibly some other edges, as was noted in \cite{di_francesco/zuber:1990} (to obtain the graph $\mathcal{E}_5^{(12)}$ from the McKay graph for the subgroup $\mathrm{K} \cong TPSL(2,7)$ an extra edge must also be added).
However, unlike with $SU(2)$, for $SU(3)$ there is a certain mismatch between the subgroups $\Gamma$, with their associated McKay graphs $\mathcal{G}_{\Gamma}$, and the $\mathcal{ADE}$ graphs. The correspondence is as in Table \ref{Table:ADE-subgroupsSU(3)}, where we use the same notation as Yau and Yu \cite{yau/yu:1993} for the subgroups E-I. The notation $\lfloor x \rfloor$ denotes the integer part of $x$. The subgroup $TA_5$ (respectively $TA_6$, $TPSL(2,7)$) is the ternary $A_5$ group (respectively ternary $A_6$ group, ternary $PSL(2,7)$ group), which is the extension of $A_5$ (respectively $A_6$, $PSL(2,7)$) by a cyclic group of order 3.
The McKay graphs are illustrated in Figures \ref{fig:extended-Agraph}-\ref{fig:extended-Kgraph}.

We will now consider the spectral measure for the McKay graph $\mathcal{G}_{\Gamma}$ associated to a finite subgroup $\Gamma \subset SU(3)$.
Any eigenvalue of $\Gamma$ can be written in the form $\chi_{\rho}(g) = \mathrm{Tr}(\rho(g))$, where $g$ is any element of the conjugacy class $\Gamma_j$ \cite{evans/pugh:2009v}.

Every element $g \in \Gamma$ is conjugate to an element $d$ in the torus, i.e. $\rho(h^{-1}gh) = \rho(d) = \mathrm{diag}(t_1,t_2,\overline{t_1 t_2})$ for some $(t_1,t_2) \in \mathbb{T}^2$. Now $\mathrm{Tr}(\rho(g)) = \mathrm{Tr}(\rho(d)) = t_1 + t_2 + \overline{t_1 t_2}$, thus its eigenvalues are all of the form $e^{i \theta_1} + e^{i \theta_2} + e^{-i (\theta_1 + \theta_2)}$, for $0 \leq \theta_1, \theta_2 < 2 \pi$, and hence the spectrum is contained in the discoid $\mathfrak{D}$.
For $\mathbb{T}^2$ or the group $SU(3)$ itself the spectrum is the whole of $\mathfrak{D}$ \cite[Section 6.2]{evans/pugh:2009v}.
So if $\Gamma$ is $SU(3)$ or one of its finite subgroups, the spectrum $\sigma(\Delta)$ of $\Delta$ is contained in $\mathfrak{D}$, illustrated in Figure \ref{fig:hypocycloid-S}.
Thus the support of $\mu_{\Delta}$ is contained in the discoid $\mathfrak{D}$.
Let $\Gamma$ be a finite subgroup of $SU(3)$ and $\Gamma_j$ its conjugacy classes, $j=1,\ldots,s$. Since the $S$-matrix simultaneously diagonalizes the representations of $\Gamma$ \cite{kawai:1989}, then as in \cite[Section 4]{evans/pugh:2009v} for the subgroups of $SU(2)$, the elements $y_i$ in (\ref{eqn:moments_general_SU(3)graph}) are then given by $y_i = S_{0,j} = \sqrt{|\Gamma_j|} \chi_0(\Gamma_j)/\sqrt{|\Gamma|} = \sqrt{|\Gamma_j|}/\sqrt{|\Gamma|}$.
Then the $m,n^{\mathrm{th}}$ moment $\varsigma_{m,n}$ is given by
\begin{equation} \label{eqn:moments-subgroupSU(3)}
\varsigma_{m,n} \; = \; \int_{\mathfrak{D}} z^m \overline{z}^n \mathrm{d}\mu(z) \; = \; \sum_{j=1}^s \frac{|\Gamma_j|}{|\Gamma|} \chi_{\rho} (\Gamma_j)^m \overline{\chi_{\rho} (\Gamma_j)}^n.
\end{equation}
Let $\Phi:\mathbb{T}^2 \rightarrow \mathfrak{D}$ be the map defined in (\ref{Phi:T^2->S}). We wish to compute `inverse' maps $\Phi^{-1}:\mathfrak{D} \rightarrow \mathbb{T}^2$ such that $\Phi \circ \Phi^{-1} = \mathrm{id}$.
For $z \in \mathfrak{D}$, we can write $z = \omega_1 + \omega_2^{-1} + \omega_1^{-1} \omega_2$ and $\overline{z} = \omega_1^{-1} + \omega_2 + \omega_1 \omega_2^{-1}$. Multiplying the first equation through by $\omega_1$, we obtain $z \omega_1 = \omega_1^2 + \omega_1 \omega_2^{-1} + \omega_2$. Then we need to find solutions $\omega_1$ to the cubic equation
\begin{equation} \label{eqn:cubic_w1}
\omega_1^3 - z \omega_1^2 + \overline{z} \omega_1 - 1 = 0.
\end{equation}
Similarly, we need to find solutions $\omega_2$ to the cubic equation $\omega_2^3 - \overline{z} \omega_2^2 + z \omega_2 - 1 = 0$. We see that the three solutions for $\omega_2$ are given by the complex conjugate of the three solutions for $\omega_1$.
Solving (\ref{eqn:cubic_w1}) we obtain solutions $\omega^{(k)}$, $k = 0,1,2$, given by \cite[Chapter V, \S6]{birkhoff/mac_lane:1965}:
$$\omega^{(k)} = (z + 2^{-1/3} \epsilon_k P + 2^{1/3} \overline{\epsilon_k} (z^2-3\overline{z}) P^{-1})/3,$$
where $\epsilon_k = e^{2 \pi i k /3}$, $2^{1/3}$ takes a real value, and $P$ is the cube root $P = (27 - 9z\overline{z} + 2z^3 + 3 \sqrt{3} \sqrt{27 - 18z\overline{z} + 4z^3 + 4\overline{z}^3 - z^2\overline{z}^2})^{1/3}$ such that $P \in \{ r e^{i \theta} | \; 0 \leq \theta < 2 \pi/3 \}$. For the roots of a cubic equation, it does not matter whether the square root in $P$ is taken to be positive or negative. We will take it to have positive value. We notice that the Jacobian $J$ appears in the expression for $P$ as the discriminant of the cubic equation (\ref{eqn:cubic_w1}).

We can define maps $\Phi^{-1}_{k,l}:\mathfrak{D} \rightarrow \mathbb{T}^2$ by
\begin{equation} \label{def:inverse_phi-SU(3)}
\Phi^{-1}_{k,l}(z) = (\omega^{(k)},\overline{\omega^{(l)}}), \qquad k, l \in \{ 0,1,2 \},
\end{equation}
for $z \in \mathfrak{D}$. Now $\Phi(\Phi^{-1}_{k,k}(z)) \neq z$, for $k =0,1,2$, however, for the other six cases ($k, l \in \{ 0,1,2 \}$ such that $k \neq l$) we do indeed have $\Phi \circ \Phi^{-1}_{k,l} = \mathrm{id}$.
These six $\Phi^{-1}_{k,l}(z)$ are the $S_3$-orbit of $\Phi^{-1}_{0,1}(z)$ under the action of the group $S_3$.
The spectral measure of $\Gamma$ (over $\mathbb{T}^2$) can then be taken as the average over these $\Phi^{-1}_{k,l}(z)$:
\begin{eqnarray}
\lefteqn{ \int_{\mathbb{T}^2} R_{m,n}(\omega_1,\omega_2) \mathrm{d}\varepsilon(\omega_1,\omega_2) } \nonumber \\
& = & \frac{1}{6} \sum_{j=1}^s \sum_{\stackrel{k,l \in \{ 0,1,2 \}:}{\scriptscriptstyle{k \neq l}}} \frac{|\Gamma_j|}{|\Gamma|} (\omega^{(k,j)} + \overline{\omega^{(l,j)}} + \overline{\omega^{(k,j)}}\omega^{(l,j)})^m (\overline{\omega^{(k,j)}} + \omega^{(l,j)} + \omega^{(k,j)}\overline{\omega^{(l,j)}})^n, \qquad \label{eqn:moments-subgroupSU(3)2}
\end{eqnarray}
where $\omega^{(p,j)}$, $j=1,\ldots,s$, are given by $\Phi_{k,l}^{-1}(\chi_{\rho}(\Gamma_j)) = (\omega^{(k,j)},\overline{\omega^{(l,j)}})$, for $k, l \in \{ 0,1,2 \}$, $k \neq l$.

\subsection{Groups A: $\mathbb{Z}_{p} \times \mathbb{Z}_{q}$}

We will now compute the spectral measure for the graph $\mathcal{G}_{\Gamma}$ corresponding to the subgroup $\Gamma = \mathbb{Z}_{p} \times \mathbb{Z}_{q}$. When $p = q = n-2$, the $SU(3)$ McKay graph of $\mathbb{Z}_{n-2} \times \mathbb{Z}_{n-2}$, Figure \ref{fig:extended-Agraph}, is the ``affine'' version of the graph $\mathcal{A}^{(n)}$ \cite[Figure 11]{behrend/pearce/petkova/zuber:2000}.
The group contains $|\Gamma| = pq$ elements, each of which is a separate conjugacy class $\Gamma_{k,l}$, where $k \in \{ 0,1,2,\ldots,p \}$, $l \in \{ 0,1,2,\ldots,q \}$. Now $\chi_{\rho}(\Gamma_{k,l}) = \widetilde{\omega}_1^{k} + \widetilde{\omega}_2^{-l} + \widetilde{\omega}_1^{-k} \widetilde{\omega}_2^{l} \in \mathfrak{D}$, where $\widetilde{\omega}_1 = e^{2\pi i/p}$, $\widetilde{\omega}_2 = e^{2\pi i/q}$. Let
\begin{equation} \label{def:Omega(k,l)}
\Omega(k,l) = (\widetilde{\omega}^{k} + \widetilde{\omega}^{-l} + \widetilde{\omega}^{l-k})^{m} (\widetilde{\omega}^{-k} + \widetilde{\omega}^{l} + \widetilde{\omega}_1^{k-l})^{n}.
\end{equation}
Then by (\ref{eqn:moments-subgroupSU(3)2}),
$$\int_{\mathbb{T}^2} R_{m,n}(\omega_1,\omega_2) \mathrm{d}\varepsilon(\omega_1,\omega_2) = \sum_{k=0}^{p-1} \sum_{l=0}^{q-1} \frac{1}{pq} \Omega(k,l),$$
and we easily obtain:

\begin{Thm}
For $\Gamma = \mathbb{Z}_{p} \times \mathbb{Z}_{q}$, the spectral measure of $\mathcal{G}_{\Gamma}$ (over $\mathbb{T}^2$) is given by the product measure
$$\mathrm{d}\varepsilon(\omega_1,\omega_2) = \mathrm{d}_{p/2} \omega_1 \; \mathrm{d}_{q/2} \omega_2,$$
where $\mathrm{d}_m$ is the uniform measure on the $2m^{\textrm{th}}$ roots of unity.
\end{Thm}

\subsection{Groups C: $\Delta(3n^2) = (\mathbb{Z}_n \times \mathbb{Z}_n) \rtimes \mathbb{Z}_3$} \label{sect:Delta(3n^2)}

The ternary trihedral group $\Delta(3n^2)$ is the semi-direct product of $\mathbb{Z}_n \times \mathbb{Z}_n$ with $\mathbb{Z}_3$, where the action of $\mathbb{Z}_3$ on $\mathbb{Z}_n \times \mathbb{Z}_n$ is given by left multiplication of $(\omega_1,\omega_2) \in \mathbb{Z}_3^2$ by the matrix $T_3$ defined in (\ref{T1,T2}), see \cite{luhn/nasri/ramond:2007}.
This group has order $|\Gamma| = 3n^2$. It has a presentation generated by the following matrices in $SU(3)$:
$$S_1 = \left( \begin{array}{ccc} 1   &   0       &   0   \\
                                0   & \omega    &   0   \\
                                0   &   0       & \omega^2   \end{array} \right),
S_2 = \left( \begin{array}{ccc} \omega   &   0       &   0   \\
                                0   & \omega    &   0   \\
                                0   &   0       & \omega   \end{array} \right)
\quad \textrm{ and } \quad
T = \left( \begin{array}{ccc} 0   &   1   &   0   \\
                                0   &   0   &   1   \\
                                1   &   0   &   0   \end{array} \right).$$
In this presentation, the action of $\mathbb{Z}_3$ on $\mathbb{Z}_n \times \mathbb{Z}_n \cong \langle S_1, S_2 \rangle$ is given by left multiplication of the element $D \in \langle S_1, S_2 \rangle$ by the matrix $T$. Here $\langle g_1, g_2, \ldots, g_s \rangle$ denotes the group generated by the elements $g_1, g_2, \ldots, g_s$.

We will consider the cases where $n \equiv 0 \textrm{ mod } 3$, $n \not \equiv 0 \textrm{ mod } 3$ separately. \\

\noindent \textbf{$n \equiv 0 \textrm{ mod } 3$} \\

\begin{figure}[htbp]
\begin{center}
  \includegraphics[width=40mm]{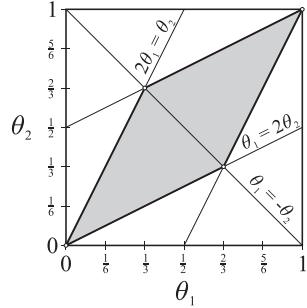}\\
 \caption{Subset $K_n \subset \mathbb{T}^2$} \label{fig:Kn}
\end{center}
\end{figure}

First we consider the case where $n \equiv 0 \textrm{ mod } 3$. The $SU(3)$ McKay graph of the group $\Delta(3n^2)$ (not drawn here) is the ``affine'' version of the graph $\mathcal{D}^{(n)}$ \cite[Figure 11]{behrend/pearce/petkova/zuber:2000}. The values of the character of the fundamental representation evaluated over the conjugacy classes of $\Delta(3n^2)$ are given in Table \ref{table:Character_table-D(3n2)} (see \cite{luhn/nasri/ramond:2007}). Here $K_n$ is the region illustrated in Figure \ref{fig:Kn} and defined by $(\theta_1,\theta_2) \in [0,1]^2 \setminus \{ (0,0), (1/3,2/3), (2/3,1/3) \}$ where $n\theta_1, n\theta_2 \in \mathbb{Z}$ and such that $2\theta_1 - \theta_2 \geq 0$, $2\theta_2 - \theta_1 \geq 0$ for $0 \leq \theta_1,\theta_2 \leq 1/2$, and $2\theta_1 - \theta_2 \leq 0$, $2\theta_2 - \theta_1 \leq 0$ for $1/2 < \theta_1,\theta_2 < 1$. Note that $|K_n| = (n^2-3)/3$
The final row in the table denotes the pair $(\theta_1, \theta_2)$ given by $(e^{2 \pi i \theta_1}, e^{2 \pi i \theta_2}) = \Phi^{-1}(\chi_{\rho}(\Gamma_j))$.

\renewcommand{\arraystretch}{1.5}

\begin{table}[bt]
\begin{center}
\begin{tabular}{|c||c|c|c|c|c|} \hline
$\Gamma_j$ & $\Gamma_1$ & $\Gamma_2$ & $\Gamma_3$ & $\Gamma_{nk,nl}$, $(k,l) \in K_n$ & $\Gamma_j'$, $j=1,\ldots,6$ \\
\hline $|\Gamma_j|$ & 1 & 1 & 1 & 3 & $n^2/3$ \\
\hline $\chi_{\rho}(\Gamma_j) \in \mathfrak{D}$ & 3 & $3 \omega$ & $3 \overline{\omega}$ & $e^{2 \pi i k} + e^{-2 \pi i l} + e^{2 \pi i (l-k)}$ & 0 \\
\hline $\Phi^{-1}(\chi_{\rho}(\Gamma_j)) \in \mathbb{T}^2$ & (1,1) & $(\omega, \overline{\omega})$ & $(\overline{\omega}, \omega)$ & $(e^{2 \pi i k}, e^{2 \pi i l})$ & $(\omega, 1)$ \\
\hline $(\theta_1, \theta_2) \in [0,1]^2$ & (0,0) & $(\frac{1}{3}, \frac{2}{3})$ & $(\frac{2}{3}, \frac{1}{3})$ & $(k,l)$ & $(\frac{1}{3},0)$ \\
\hline
\end{tabular} \\
\caption{$\chi_{\rho}(\Gamma_j)$ for group $\Delta(3n^2)$, $n \equiv 0 \textrm{ mod } 3$. Here $\omega = e^{2 \pi i/3}$.} \label{table:Character_table-D(3n2)}
\end{center}
\end{table}

\renewcommand{\arraystretch}{1}

Let $\widetilde{\omega} = e^{2 \pi i/n}$, and $\Omega(k,l)$ be as in (\ref{def:Omega(k,l)}).
Then by (\ref{eqn:moments-subgroupSU(3)2}),
\begin{eqnarray}
\lefteqn{ \int_{\mathbb{T}^2} R_{m_1,m_2}(\omega_1,\omega_2) \mathrm{d}\varepsilon(\omega_1,\omega_2) } \nonumber \\
& = & \frac{1}{3n^2} \Omega(0,0) + \frac{1}{3n^2} \Omega(\textstyle \frac{1}{3}, \frac{2}{3} \displaystyle) + \frac{1}{3n^2} \Omega(\textstyle \frac{2}{3}, \frac{1}{3} \displaystyle) + \frac{3}{3n^2} \sum_{k,l \in K_n} \Omega(k,l) + \frac{n^2/3}{3n^2} \sum_{j=1}^{6} \Omega(\textstyle \frac{1}{3}, 0 \displaystyle) \nonumber \\
& = & \frac{1}{3n^2} \Omega(0,0) + \frac{1}{3n^2} \Omega(\textstyle \frac{1}{3}, \frac{2}{3} \displaystyle) + \frac{1}{3n^2} \Omega(\textstyle \frac{2}{3}, \frac{1}{3} \displaystyle) + \frac{1}{3n^2} \sum_{k,l} \Omega(k,l) + \frac{1}{9} \sum_{g \in S_3} \Omega(g(\textstyle \frac{1}{3} \displaystyle,0)) \label{eqn:Delta3n^2-I} \\
& = & \frac{1}{3n^2} \sum_{nk, nl \in \mathbb{Z}_n} \Omega(k,l) + \frac{1}{9} \sum_{g \in S_3} \Omega(g(\textstyle \frac{1}{3} \displaystyle,0)), \nonumber
\end{eqnarray}
where the first summation in (\ref{eqn:Delta3n^2-I}) is over all $(k,l) \neq (0,0), (1/3,2/3), (2/3,1/3)$ such that $nk, nl \in \mathbb{Z}_n$.
Now
$$\sum_{g \in S_3} \Omega(g(0, \textstyle \frac{1}{3} \displaystyle)) = \frac{1}{4 \pi^4} \int_{\mathbb{T}^2} R_{m_1,m_2}(\omega_1,\omega_2) J^2 \mathrm{d}^{(3)}(\omega_1,\omega_2),$$
since $J(\omega_1,\omega_2) = 0$ for $(\omega_1,\omega_2) \in \mathbb{T}^2$ such that $\Phi(\omega_1,\omega_2)$ is on the deltoid, the boundary of the discoid $\mathfrak{D}$ (c.f. Section \ref{sec:E4(12)}).
Then the spectral measure (over $\mathbb{T}^2$) for the group $\Delta(3n^2)$, $n \equiv 0 \textrm{ mod } 3$, is
$$\mathrm{d}\varepsilon(\omega_1,\omega_2) = \frac{1}{3} \, \mathrm{d}_{n/2}\omega_1 \; \mathrm{d}_{n/2}\omega_2 + \frac{1}{36 \pi^4} J^2 \mathrm{d}^{(3)}(\omega_1,\omega_2),$$
where $\mathrm{d}_n$ is the uniform measure over $2n^{\mathrm{th}}$ roots of unity, $\omega = e^{2 \pi i/3}$, and $\mathrm{d}^{(n)}$ is the uniform measure over the points in $D_n$. \\

\noindent \textbf{$n \not \equiv 0 \textrm{ mod } 3$} \\

\renewcommand{\arraystretch}{1.5}

\begin{table}[bt]
\begin{center}
\begin{tabular}{|c||c|c|c|} \hline
$\Gamma_j$ & $\Gamma_1$ & $\Gamma_{nk,nl}$, $(k,l) \in K_n$ & $\Gamma_j'$, $j=1,2$ \\
\hline $|\Gamma_j|$ & 1 & 3 & $n^2$ \\
\hline $\chi_{\rho}(\Gamma_j) \in \mathfrak{D}$ & 3 & $e^{2 \pi i k} + e^{-2 \pi i l} + e^{2 \pi i (l-k)}$ & 0  \\
\hline $\Phi^{-1}(\chi_{\rho}(\Gamma_j)) \in \mathbb{T}^2$ & (1,1) & $(e^{2 \pi i k}, e^{2 \pi i l})$ & $(\omega, 1)$ \\
\hline $(\theta_1, \theta_2) \in [0,1]^2$ & (0,0) & $(k,l)$ & $(\frac{1}{3},0)$ \\
\hline
\end{tabular} \\
\caption{$\chi_{\rho}(\Gamma_j)$ for group $\Delta(3n^2)$, $n \not \equiv 0 \textrm{ mod } 3$.} \label{table:Character_table-D(3n2)-2}
\end{center}
\end{table}

\renewcommand{\arraystretch}{1}

Now consider the case where $n \not \equiv 0 \textrm{ mod } 3$. The values of $\chi_{\rho}(\Gamma_j)$ for $\Delta(3n^2)$ are given in Table \ref{table:Character_table-D(3n2)-2} (see \cite{luhn/nasri/ramond:2007}). The final row in the table denotes the pair $(\theta_1, \theta_2)$ given by $(e^{2 \pi i \theta_1}, e^{2 \pi i \theta_2}) = \Phi^{-1}(\chi_{\rho}(\Gamma_j))$.
Then by (\ref{eqn:moments-subgroupSU(3)2}),
\begin{eqnarray*}
\lefteqn{ \int_{\mathbb{T}^2} R_{m_1,m_2}(\omega_1,\omega_2) \mathrm{d}\varepsilon(\omega_1,\omega_2) } \nonumber \\
& = & \frac{1}{3n^2} \Omega(0,0) + \frac{3}{3n^2} \sum_{k,l \in K_n} \Omega(k,l) + \frac{n^2}{3n^2} \sum_{j=1}^{2} \Omega(\textstyle \frac{1}{3}, 0 \displaystyle) \; = \; \frac{1}{3n^2} \sum_{nk, nl \in \mathbb{Z}_n} \Omega(k,l) + \frac{2}{3} \Omega(\textstyle \frac{1}{3}, 0 \displaystyle) \\
& = & \frac{1}{3n^2} \sum_{nk, nl \in \mathbb{Z}_n} \Omega(k,l) + \frac{1}{9} \sum_{g \in S_3} \Omega(g(\textstyle \frac{1}{3} \displaystyle,0)),
\end{eqnarray*}
and we obtain the same spectral measure as for the case $n \equiv 0 \textrm{ mod } 3$. Summarizing, we have:

\begin{Thm}
The spectral measure (over $\mathbb{T}^2$) for the ternary trihedral group $\Delta(3n^2) = (\mathbb{Z}_n \times \mathbb{Z}_n) \rtimes \mathbb{Z}_3$ is
\begin{equation}
\mathrm{d}\varepsilon(\omega_1,\omega_2) = \frac{1}{3} \, \mathrm{d}_{n/2}\omega_1 \; \mathrm{d}_{n/2}\omega_2 + \frac{1}{36 \pi^4} J^2 \mathrm{d}^{(3)}(\omega_1,\omega_2),
\end{equation}
where $\mathrm{d}_n$ is the uniform measure over $2n^{\mathrm{th}}$ roots of unity and $\mathrm{d}^{(3)}$ is the uniform measure over the points in $D_3$.
\end{Thm}

\noindent
\emph{Remark:}
The spectral measure of the binary dihedral group $BD_n \subset SU(2)$ is \cite[Theorem 4.1]{banica/bisch:2007}, \cite[Section 4.2]{evans/pugh:2009v}:
$$\mathrm{d} \varepsilon(u) = \frac{1}{2} \mathrm{d}_{n-2} u + \frac{1}{4}(\delta_i + \delta_{-i}),$$
where $\delta_x$ is the Dirac measure at the point $x$. The points $i$, $-i$ are the two points in $\mathbb{T}$ which map to zero in the interval $[-2,2]$ under the map $\Phi$ from $\mathbb{T}$ to the support of the spectral measure for (a subgroup of) $SU(2)$. See \cite{evans/pugh:2009v} for more details.
Let $\omega = e^{2 \pi i/3}$ as before. Since
$$\sum_{g \in S_3} \Omega(g(\textstyle \frac{1}{3} \displaystyle,0)) = \Omega(0, \textstyle \frac{1}{3} \displaystyle) + \Omega(0, \textstyle \frac{2}{3} \displaystyle) + \Omega(\textstyle \frac{1}{3} \displaystyle,0) + \Omega(\textstyle \frac{1}{3}, \frac{1}{3} \displaystyle) + \Omega(\textstyle \frac{2}{3} \displaystyle,0) + \Omega(\textstyle \frac{2}{3}, \frac{2}{3} \displaystyle)$$
it is easy to see that
$$\mathrm{d}\varepsilon(\omega_1,\omega_2) = \frac{1}{3} \, \mathrm{d}_{n/2}\omega_1 \; \mathrm{d}_{n/2}\omega_2 + \frac{\delta_{(1,\omega)} + \delta_{(1,\overline{\omega})} + \delta_{(\omega,1)} + \delta_{(\omega,\omega)} + \delta_{(\overline{\omega},1)} + \delta_{(\overline{\omega},\overline{\omega})}}{9},$$
where now the points $(1,\omega)$, $(1,\overline{\omega})$, $(\omega,1)$, $(\omega,\omega)$, $(\overline{\omega},1)$, $(\overline{\omega},\overline{\omega})$ are the points in $\mathbb{T}^2$ which map to zero in the support $\mathfrak{D}$ of the spectral measure for (a subgroup of) $SU(3)$ under the map $\Phi$ defined in (\ref{Phi:T^2->S}).

\subsection{Groups D: $\Delta(6n^2) = (\mathbb{Z}_n \times \mathbb{Z}_n) \rtimes S_3$}

The ternary trihedral group $\Delta(6n^2)$ is the semi-direct product of $\mathbb{Z}_n \times \mathbb{Z}_n$ with $S_3$, where the action of $S_3$ on $\mathbb{Z}_n \times \mathbb{Z}_n$ is given by left multiplication of $(\omega_1,\omega_2) \in \mathbb{Z}_3^2$ by the matrices $T_2$, $T_3$ defined in (\ref{T1,T2}) which generate $S_3$, see \cite{escobar/luhn:2009}.
This group has order $|\Gamma| = 6n^2$. It has a presentation generated by the generators $S_1, S_2$ and $T$ of $\Delta(3n^2)$, and the matrix $Q \in SU(3)$ given by
$$Q = \left( \begin{array}{ccc}   -1  &   0   &   0   \\
                                  0   &   0   &   -1  \\
                                  0   &   -1  &   0   \end{array} \right).$$
In this presentation, the action of $S_3$ on $\mathbb{Z}_n \times \mathbb{Z}_n \cong \langle S_1, S_2 \rangle$ is given by left multiplication of the element $D \in \langle S_1, S_2 \rangle$ by the matrices $T, Q$.

We will consider the cases where $n \equiv 0 \textrm{ mod } 3$, $n \not \equiv 0 \textrm{ mod } 3$ separately. \\

\noindent \textbf{$n \equiv 0 \textrm{ mod } 3$} \\

First we consider the case where $n \equiv 0 \textrm{ mod } 3$.
The values of $\chi_{\rho}(\Gamma_j)$ for $\Delta(6n^2)$ is given in Table \ref{table:Character_table-D(6n2)} (see \cite{escobar/luhn:2009}), where $\theta_1^{(k)}$, $\theta_2^{(k)}$ are defined by
\begin{equation} \label{eqn:theta(k)}
(\theta_1^{(k)}, \theta_1^{(k)}) = \left\{ \begin{array}{cl}
                                            (\frac{1}{2}k + \frac{1}{4}, k) & \textrm{ if } \frac{1}{6} \leq k < \frac{1}{2}, \\
                                            (k, \frac{1}{2}k + \frac{1}{4}) & \textrm{ if } \frac{1}{2} \leq k < \frac{5}{6}, \\
                                            (\frac{1}{2}k + \frac{1}{4}, \frac{1}{4} - \frac{1}{2}k) & \textrm{ otherwise.} \\
                                            \end{array} \right.
\end{equation}
The set $K_n'$ is the subset of $K_n$ given by $(\theta_1,\theta_2) \in K_n$ such that $\theta_1 + \theta_2 < 1$, $2\theta_1 - \theta_2 < 0$ and $2\theta_2 - \theta_1 < 0$.

\renewcommand{\arraystretch}{1.5}

\begin{table}[bt]
\begin{center}
\begin{tabular}{|c||c|c|c|c|} \hline
$\Gamma_j$ & $\Gamma_1$ & $\Gamma_2$ & $\Gamma_3$ & $\Gamma_{(nk)}$, $nk \in \mathbb{Z}_n \setminus \{ 0, \textstyle \frac{n}{3} \displaystyle, \textstyle \frac{2n}{3} \displaystyle \}$ \\
\hline $|\Gamma_j|$ & 1 & 1 & 1 & 3 \\
\hline $\chi_{\rho}(\Gamma_j) \in \mathfrak{D}$ & 3 & $3 \omega$ & $3 \overline{\omega}$ & $e^{-4 \pi i k} + 2e^{2 \pi i k}$ \\
\hline $(e^{2 \pi i \theta_1}, e^{2 \pi i \theta_2}) = \Phi^{-1}(\chi_{\rho}(\Gamma_j)) \in \mathbb{T}^2$ & (1,1) & $(\omega, \overline{\omega})$ & $(\overline{\omega}, \omega)$ & $(e^{2 \pi i k}, e^{-2 \pi i k})$ \\
\hline $(\theta_1, \theta_2) \in [0,1]^2$ & (0,0) & $(\frac{1}{3}, \frac{2}{3})$ & $(\frac{2}{3}, \frac{1}{3})$ & $(k,-k)$ \\
\hline
\end{tabular}
$$\;$$
\begin{tabular}{|c||c|c|c|} \hline
$\Gamma_j$ & $\Gamma_{nk,nl}$, $(k,l) \in K_n' \setminus \{ (\textstyle \frac{1}{3}, \frac{1}{3} \displaystyle) \}$ & $\Gamma_j'$, $j=1,2,3$ & $\Gamma_{(nk)}'$, $nk \in \mathbb{Z}_n$ \\
\hline $|\Gamma_j|$ & 6 & $2n^2/3$ & $3n$ \\
\hline $\chi_{\rho}(\Gamma_j) \in \mathfrak{D}$ & $e^{2 \pi i k} + e^{-2 \pi i l} + e^{2 \pi i (l-k)}$ & 0 & $e^{-2 \pi i k}$ \\
\hline $\Phi^{-1}(\chi_{\rho}(\Gamma_j)) \in \mathbb{T}^2$ & $(e^{2 \pi i k}, e^{2 \pi i l})$ & $(\omega, 1)$ & $(e^{2 \pi i \theta_1^{(k)}}, e^{2 \pi i \theta_2^{(k)}})$ \\
\hline $(\theta_1, \theta_2) \in [0,1]^2$ & $(k,l)$ & $(\frac{1}{3},0)$ & $(\theta_1^{(k)}, \theta_2^{(k)})$ \\
\hline
\end{tabular} \\
\caption{$\chi_{\rho}(\Gamma_j)$ for group $\Delta(6n^2)$, $n \equiv 0 \textrm{ mod } 3$. Here $\omega = e^{2 \pi i/3}$.} \label{table:Character_table-D(6n2)}
\end{center}
\end{table}

\renewcommand{\arraystretch}{1}

\begin{figure}[bt]
\begin{center}
  \includegraphics[width=105mm]{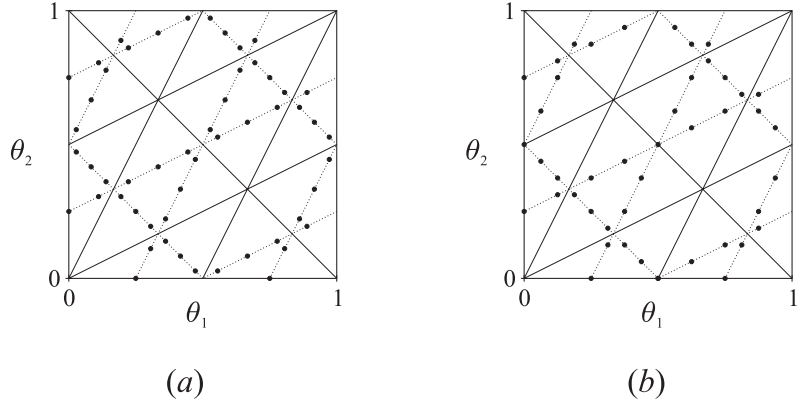}\\
 \caption{The points $(\theta_1^{(k)},\theta_2^{(k)})$ for $(a)$ $n = 9$, $(b)$ $n = 8$.} \label{fig:poly-17}
\end{center}
\end{figure}

The points $(\theta_1, \theta_2) = (k,-k)$, for all $k \neq 1/3, 2/3$ such that $nk \in \mathbb{Z}_n$, lie on the boundary of the fundamental domain $C$ in Figure \ref{fig:poly-4} and have weight $|\Gamma_j|/|\Gamma| = 3/6n^2$, thus the $S_3$-orbit of $(k,-k)$ contains the three points $(k,-k)$, $(k,2k)$ and $(-2k,-k)$, each of which has weight $1/6n^2$. The points $(k,l) \in K_n'$ have weight $|\Gamma_j|/|\Gamma| = 6/6n^2$, thus the six points in the $S_3$-orbit of $(k,l) \in K_n'$ each have weight $1/6n^2$ also. Furthermore, the $S_3$-orbit of $(1/3, 0)$ contains six points.
Then we see from (\ref{eqn:moments-subgroupSU(3)2}) that
\begin{eqnarray*}
\lefteqn{ \int_{\mathbb{T}^2} R_{m_1,m_2}(\omega_1,\omega_2) \mathrm{d}\varepsilon(\omega_1,\omega_2) } \\
& = & \frac{1}{6n^2} \Omega(0,0) + \frac{1}{6n^2} \Omega(\textstyle \frac{1}{3}, \frac{2}{3} \displaystyle) + \frac{1}{6n^2} \Omega(\textstyle \frac{2}{3}, \frac{1}{3} \displaystyle) + \frac{1}{6n^2} \sum_{k,l} \Omega(k,l) + \frac{2n^2/3}{6n^2} \sum_{j=1}^{3} \Omega(\textstyle \frac{1}{3}, 0 \displaystyle) \\
& & + \frac{3n}{6n^2} \sum_{nk \in \mathbb{Z}_n} \Omega(\theta_1^{(k)}, \theta_2^{(k)}) \\
& = & \frac{1}{6n^2} \sum_{nk, nl \in \mathbb{Z}_n} \Omega(k,l) + \frac{1}{18} \sum_{g \in S_3} \Omega(g(\textstyle \frac{1}{3} \displaystyle,0)) + \frac{1}{12n} \sum_{g \in S_3} \sum_{nk \in \mathbb{Z}_n} \Omega(g(\theta_1^{(k)}, \theta_2^{(k)})),
\end{eqnarray*}
where in the first equality the first summation is over all $k,l$ such that $nk,nl \in \mathbb{Z}_n$ and $(k,l) \neq (0,0), (1/3,2/3), (2/3,1/3)$.
In the second equality, the second summation $\sum_{g \in S_3} \Omega(g(\textstyle \frac{1}{3} \displaystyle,0))$ is given by the measure $J^2 \mathrm{d}^{(3)}$ as in Section \ref{sect:Delta(3n^2)}.
The points $(e^{2 \pi i \theta_1^{(k)}}, e^{2 \pi i \theta_2^{(k)}})$ given by the third summation $\sum_{g \in S_3} \sum_{nk \in \mathbb{Z}_n} \Omega(g(\theta_1^{(k)}, \theta_2^{(k)}))$ lie on the lines $\theta_1 + \theta_2 = 1/2$, $2\theta_1 - \theta_2 = 1/2$ and $2\theta_2 - \theta_1 = 1/2$, where $n$ points lie equidistantly along the length of each line such that there is a point at each of $(1/4,1/4)$, $(1/4,0)$ and $(0,1/4)$. These points are illustrated in Figure \ref{fig:poly-17}$(a)$ for $n = 9$. There are $6n$ distinct points when $n \not \equiv 0 \textrm{ mod } 6$, but only $6n-9$ distinct points when $n \equiv 0 \textrm{ mod } 6$ as the points $(-1,-1)$, $(e^{\pi i/3},\omega)$, $(\omega, e^{\pi i/3})$, and their $S_3$-orbits, have multiplicity two.
The third summation is thus given by a different measure $\mathrm{d}\varepsilon_n'$ depending on whether $n$ is divisible by 6 or not. When $n \equiv 0 \textrm{ mod } 6$,
$$\mathrm{d}\varepsilon_n' = 18 \mathrm{d}^{((4))} + 36 \sum_{j=1}^{\lfloor (n+3)/6 \rfloor} \mathrm{d}^{(4n/(n-2j),j/n)},$$
whilst for $n \equiv 0 \textrm{ mod } 3$, $n \not \equiv 0 \textrm{ mod } 6$,
$$\mathrm{d}\varepsilon_n' = 18 \mathrm{d}^{((4))} + 36 \sum_{j=1}^{\lfloor (n+3)/6 \rfloor} \mathrm{d}^{(4n/(n-2j),j/n)} + 18\mathrm{d}^{((2))}.$$

Thus, we have:

\begin{Thm}
The spectral measure (over $\mathbb{T}^2$) for the ternary trihedral group $\Delta(6n^2) = (\mathbb{Z}_n \times \mathbb{Z}_n) \rtimes S_3$, $n \equiv 0 \textrm{ mod } 3$, is
\begin{equation}
\mathrm{d}\varepsilon = \frac{1}{6} \, \mathrm{d}_{n/2} \; \mathrm{d}_{n/2} + \frac{1}{72 \pi^4} J^2 \, \mathrm{d}^{(3)} + \frac{3}{2n} \mathrm{d}^{((4))} + \frac{3}{n} \sum_{j=1}^{\lfloor (n+3)/6 \rfloor} \mathrm{d}^{(4n/(n-2j),j/n)} + \mathrm{d}',
\end{equation}
where $\mathrm{d}'$ is given by
$$\mathrm{d}' = \left\{ \begin{array}{cl}
                        \frac{3}{2n} \mathrm{d}^{((2))} & \textrm{ if } n \equiv 0 \textrm{ mod } 6, \\
                        0 & \textrm{ if } n \not \equiv 0 \textrm{ mod } 6. \end{array} \right.$$
\end{Thm}

\noindent \textbf{$n \not \equiv 0 \textrm{ mod } 3$} \\

Now consider the case where $n \not \equiv 0 \textrm{ mod } 3$. The values of $\chi_{\rho}(\Gamma_j)$ for $\Delta(6n^2)$ are given in Table \ref{table:Character_table-D(6n2)-2} (see \cite{escobar/luhn:2009}). The final row in the table denotes the pair $(\theta_1, \theta_2)$ given by $(e^{2 \pi i \theta_1}, e^{2 \pi i \theta_2}) = \Phi^{-1}(\chi_{\rho}(\Gamma_j)) \in \mathbb{T}^2$, where $\theta_1^{(k)}$, $\theta_2^{(k)}$ are defined by (\ref{eqn:theta(k)}).

\renewcommand{\arraystretch}{1.5}

\begin{table}[bt]
\begin{center}
\begin{tabular}{|c||c|c|} \hline
$\Gamma_j$ & $\Gamma_1$ & $\Gamma_{(nk)}$, $nk \in \mathbb{Z}_n \setminus \{ 0 \}$ \\
\hline $|\Gamma_j|$ & 1 & 3 \\
\hline $\chi_{\rho}(\Gamma_j) \in \mathfrak{D}$ & 3 & $e^{-4 \pi i k} + 2e^{2 \pi i k}$ \\
\hline $(e^{2 \pi i \theta_1}, e^{2 \pi i \theta_2}) = \Phi^{-1}(\chi_{\rho}(\Gamma_j)) \in \mathbb{T}^2$ & (1,1) & $(e^{2 \pi i k}, e^{-2 \pi i k})$ \\
\hline $(\theta_1, \theta_2) \in [0,1]^2$ & (0,0) & $(k,-k)$ \\
\hline
\end{tabular}
$$\;$$
\begin{tabular}{|c||c|c|c|} \hline
$\Gamma_j$ & $\Gamma_{nk,nl}$, $(k,l) \in K_n'$ & $\Gamma_j'$, $j=1,2,3$ & $\Gamma_{(nk)}'$, $nk \in \mathbb{Z}_n$ \\
\hline $|\Gamma_j|$ & 6 & $2n^2$ & $3n$ \\
\hline $\chi_{\rho}(\Gamma_j) \in \mathfrak{D}$ & $e^{2 \pi i k} + e^{-2 \pi i l} + e^{2 \pi i (l-k)}$ & 0 & $e^{-2 \pi i k}$ \\
\hline $\Phi^{-1}(\chi_{\rho}(\Gamma_j)) \in \mathbb{T}^2$ & $(e^{2 \pi i k}, e^{2 \pi i l})$ & $(\omega, 1)$ & $(e^{2 \pi i \theta_1^{(k)}}, e^{2 \pi i \theta_2^{(k)}})$ \\
\hline $(\theta_1, \theta_2) \in [0,1]^2$ & $(k,l)$ & $(\frac{1}{3},0)$ & $(\theta_1^{(k)}, \theta_2^{(k)})$ \\
\hline
\end{tabular} \\
\caption{$\chi_{\rho}(\Gamma_j)$ for group $\Delta(6n^2)$, $n \not \equiv 0 \textrm{ mod } 3$. Here $\omega = e^{2 \pi i/3}$.} \label{table:Character_table-D(6n2)-2}
\end{center}
\end{table}

\renewcommand{\arraystretch}{1}

Then from (\ref{eqn:moments-subgroupSU(3)2}) we have
\begin{eqnarray*}
\lefteqn{ \int_{\mathbb{T}^2} R_{m_1,m_2}(\omega_1,\omega_2) \mathrm{d}\varepsilon(\omega_1,\omega_2) } \\
& = & \frac{1}{6n^2} \Omega(0,0) + \frac{1}{6n^2} \sum_{k,l} \Omega(k,l) + \frac{2n^2}{6n^2} \Omega(\textstyle \frac{1}{3}, 0 \displaystyle) + \frac{3n}{6n^2} \sum_{nk \in \mathbb{Z}_n} \Omega(\theta_1^{(k)}, \theta_2^{(k)}) \\
& = & \frac{1}{6n^2} \sum_{nk, nl \in \mathbb{Z}_n} \Omega(k,l) + \frac{1}{18} \sum_{g \in S_3} \Omega(g(\textstyle \frac{1}{3} \displaystyle,0)) + \frac{1}{12n} \sum_{g \in S_3} \sum_{nk \in \mathbb{Z}_n} \Omega(g(\theta_1^{(k)}, \theta_2^{(k)})),
\end{eqnarray*}
where in the first equality the first summation is over all $k,l$ such that $nk,nl \in \mathbb{Z}_n$ and $(k,l) \neq (0,0)$.
In the second equality, the points $(e^{2 \pi i \theta_1^{(k)}}, e^{2 \pi i \theta_2^{(k)}})$ given by the third summation $\sum_{g \in S_3} \sum_{nk \in \mathbb{Z}_n} \Omega(g(\theta_1^{(k)}, \theta_2^{(k)}))$ again lie on the lines $\theta_1 + \theta_2 = 1/2$, $2\theta_1 - \theta_2 = 1/2$ and $2\theta_2 - \theta_1 = 1/2$, where $n$ points lie equidistantly along the length of each line such that there is a point at each of $(1/4,1/4)$, $(1/4,0)$ and $(0,1/4)$. These points are illustrated in Figure \ref{fig:poly-17}$(b)$ for $n = 8$. There are $6n-3$ distinct points, the points $(-1,-1)$, $(-1,1)$ and $(1,-1)$ having multiplicity two. The measure $\mathrm{d}\varepsilon'$ given by the third summation cannot be written as a linear combination of the measures in Definition \ref{def:4measures}, since all the measures in Definition \ref{def:4measures} are topologically invariant under a rotation of each triangular fundamental domain of $\mathbb{T}^2/S_3$ (see Figure \ref{fig:poly-4}) by $2\pi/3$, but $\mathrm{d}\varepsilon'$ is not. The measure $\mathrm{d}\varepsilon'$ is given by
$$\mathrm{d}\varepsilon' = \sum_{j=1}^n \left( \delta_{(e^{2 \pi i j/n}, e^{2 \pi i (1 + 2j)/2n})} + \delta_{(e^{2 \pi i (1 + 2j)/2n}, e^{2 \pi i j/n})} + \delta_{(e^{2 \pi i (1 + 2j)/2n}, e^{2 \pi i (1 - 2j)/2n})} \right),$$
where $\delta_x$ is the Dirac measure at the point $x$.

Thus, we have:

\begin{Thm}
The spectral measure (over $\mathbb{T}^2$) for the ternary trihedral group $\Delta(6n^2) = (\mathbb{Z}_n \times \mathbb{Z}_n) \rtimes S_3$, $n \not \equiv 0 \textrm{ mod } 3$, is
\begin{eqnarray}
\mathrm{d}\varepsilon & = & \frac{1}{6} \, \mathrm{d}_{n/2} \; \mathrm{d}_{n/2} + \frac{1}{72 \pi^4} J^2 \, \mathrm{d}^{(3)} + \frac{3}{2n} \mathrm{d}^{((4))} + \frac{3}{n} \sum_{j=1}^{\lfloor (n+3)/6 \rfloor} \mathrm{d}^{(4n/(n-2j),j/n)} \nonumber \\
& & + \frac{1}{12n} \sum_{j=1}^n \left( \delta_{(e^{2 \pi i j/n}, e^{2 \pi i (1 + 2j)/2n})} + \delta_{(e^{2 \pi i (1 + 2j)/2n}, e^{2 \pi i j/n})} + \delta_{(e^{2 \pi i (1 + 2j)/2n}, e^{2 \pi i (1 - 2j)/2n})} \right). \quad \;
\end{eqnarray}
\end{Thm}

\subsection{Group E $= \Sigma(36 \times 3) = \Delta (3.3^2) \rtimes \mathbb{Z}_4$} \label{sect:(E)}

\renewcommand{\arraystretch}{1.5}

\begin{table}[bt]
\begin{center}
\small
\begin{tabular}{|c|ccc|ccc|} \hline
$j$ & 1 & 2 & 3 & 4, 5 & 6, 7 & 8, 9 \\
\hline $|\Gamma_j|$ & 1 & 1 & 1 & 9 & 9 & 9 \\
\hline $\chi_{\rho}(\Gamma_j) \in \mathfrak{D}$ & 3 & $3\omega$ & $3\overline{\omega}$ & 1 & $\omega$ & $\overline{\omega}$ \\
\hline $(e^{2 \pi i \theta_1}, e^{2 \pi i \theta_2}) = \Phi^{-1}(\chi_{\rho}(\Gamma_j))$ & (1,1) & $(\omega, \overline{\omega})$ & $(\overline{\omega}, \omega)$ & $(1,i)$ & $(-\omega i,-\overline{\omega} i)$ & $(-\overline{\omega} i,-\omega i)$ \\
\hline $(\theta_1, \theta_2) \in [0,1]^2$ & (0,0) & $(\frac{1}{3}, \frac{2}{3})$ & $(\frac{2}{3}, \frac{1}{3})$ & $(0, \frac{1}{4})$ & $(\frac{1}{12}, \frac{5}{12})$ & $(\frac{5}{12}, \frac{1}{12})$ \\
\hline
\end{tabular}
$$\;$$
\begin{tabular}{|c|ccc|c|} \hline
$j$ & 10 & 11 & 12 & 13, 14 \\
\hline $|\Gamma_j|$ & 9 & 9 & 9 & 12 \\
\hline $\chi_{\rho}(\Gamma_j) \in \mathfrak{D}$ & $-1$ & $-\omega$ & $-\overline{\omega}$ & 0 \\
\hline $(e^{2 \pi i \theta_1}, e^{2 \pi i \theta_2}) \in \mathbb{T}^2$ & (1,-1) & $(\omega,-\overline{\omega})$ & $(-\overline{\omega},\omega)$ & $(1,\omega)$ \\
\hline $(\theta_1, \theta_2) \in [0,1]^2$ & $(0, \frac{1}{2})$ & $(\frac{2}{6}, \frac{1}{6})$ & $(\frac{1}{6}, \frac{2}{6})$ & $(\frac{1}{3}, 0)$ \\
\hline
\end{tabular} \\
\normalsize
\caption{$\chi_{\rho}(\Gamma_j)$ for group E $= \Sigma(36 \times 3)$. Here $\omega = e^{2 \pi i/3}$.} \label{table:Character_table-(E)}
\end{center}
\end{table}

\renewcommand{\arraystretch}{1}

The subgroup E has order 108, and its McKay graph, Figure \ref{fig:extended-E(8)graph}, is the ``affine'' version of the graph $\mathcal{E}^{(8)}$ \cite[Figure 13]{evans/pugh:2009i}. It is the semidirect product of the ternary trihedral group $\Delta (3.3^2)$ with $\mathbb{Z}_4$, where the action of $\mathbb{Z}_4$ on $\Delta (3.3^2) = \langle S_1, S_2, T \rangle$ is given by left multiplication of the element $D \in \Delta (3.3^2)$ by the matrix $V \in SU(3)$ given by
$$V = \frac{1}{\sqrt{-3}} \left( \begin{array}{ccc}
                                1   &   1      &   1      \\
                                1   & \omega   & \omega^2 \\
                                1   & \omega^2 & \omega   \end{array} \right).$$
Note that $V^2 = Q$, so that E also contains the ternary trihedral group $\Delta(6.3^2)$ as a normal subgroup.
This group is a subgroup of $SU(3)$ but not of $SU(3)/\mathbb{Z}_3$
The values of $\chi_{\rho}(\Gamma_j)$ for E are given in Table \ref{table:Character_table-(E)} (see \cite{hanany/he:1999}).
Then, by (\ref{eqn:moments-subgroupSU(3)2}),
\begin{eqnarray*}
\lefteqn{ \int_{\mathbb{T}^2} R_{m_1,m_2}(\omega_1,\omega_2) \mathrm{d}\varepsilon(\omega_1,\omega_2) } \nonumber \\
& = & \frac{1}{108} \Omega(0,0) + \frac{1}{108} \Omega(\textstyle \frac{1}{3}, \frac{2}{3} \displaystyle) + \frac{1}{108} \Omega(\textstyle \frac{2}{3}, \frac{1}{3} \displaystyle) + \frac{9+9}{108} \left( \Omega(0, \textstyle \frac{1}{4} \displaystyle) + \Omega(\textstyle \frac{1}{12}, \frac{5}{12} \displaystyle) + \Omega(\textstyle \frac{5}{12}, \frac{1}{12} \displaystyle) \right) \\
& & + \frac{9}{108} \left( \Omega(0, \textstyle \frac{1}{2} \displaystyle) + \Omega(\textstyle \frac{2}{6}, \frac{1}{6} \displaystyle) + \Omega(\textstyle \frac{1}{6}, \frac{2}{6} \displaystyle) \right) + \frac{12+12}{108} \Omega(\textstyle \frac{1}{3}, 0 \displaystyle) \\
& = & \frac{1}{108} \Omega(0,0) + \frac{1}{108} \Omega(\textstyle \frac{1}{3}, \frac{2}{3} \displaystyle) + \frac{1}{108} \Omega(\textstyle \frac{2}{3}, \frac{1}{3} \displaystyle) \\
& & + \frac{1}{36} \sum_{g \in S_3} \left( \Omega(g(0, \textstyle \frac{1}{4} \displaystyle)) + \Omega(g(\textstyle \frac{1}{12}, \frac{5}{12} \displaystyle)) + \Omega(g(\textstyle \frac{5}{12}, \frac{1}{12} \displaystyle)) \right) \\
& & + \frac{1}{36} \sum_{g \in S_3} \left( \Omega(g(0, \textstyle \frac{1}{2} \displaystyle)) + \Omega(g(\textstyle \frac{2}{6}, \frac{1}{6} \displaystyle)) + \Omega(g(\textstyle \frac{1}{6}, \frac{2}{6} \displaystyle)) \right) + \frac{1}{27} \sum_{g \in S_3} \Omega(g(\textstyle \frac{1}{3}, 0 \displaystyle)).
\end{eqnarray*}
The set of the three fixed points $(0,0)$, $(1/3,2/3)$ and $(2/3,1/3)$ is $D_1$, whilst the last summation above is given by the measure $J^2 \mathrm{d}^{(3)}(\omega_1,\omega_2)/4\pi^4$ as in Section \ref{sect:Delta(3n^2)}.
We also have
\begin{eqnarray*}
\lefteqn{ \sum_{g \in S_3} \left( \Omega(g(0,1/4)) + \Omega(g(1/12,5/12)) + \Omega(g(5/12,1/12)) \right) } \\
& & \hspace{2cm} = 18 \int_{\mathbb{T}^2} R_{m_1,m_2}(\omega_1,\omega_2) \mathrm{d}^{((4))}(\omega_1,\omega_2), \hspace{4.5cm} \left. \right.
\end{eqnarray*}
and
\begin{eqnarray*}
\lefteqn{ \sum_{g \in S_3} \left( \Omega(g(0,1/2)) + \Omega(g(2/6,1/6)) + \Omega(g(1/6,2/6)) \right) } \\
& & \hspace{1cm} = 12 \int_{\mathbb{T}^2} R_{m_1,m_2}(\omega_1,\omega_2) \mathrm{d}^{(2)}(\omega_1,\omega_2) - 3 \int_{\mathbb{T}^2} R_{m_1,m_2}(\omega_1,\omega_2) \mathrm{d}^{(1)}(\omega_1,\omega_2) .
\end{eqnarray*}

Then using (\ref{eqn:measure-relation2}) we obtain the following:

\begin{Thm}
The spectral measure (over $\mathbb{T}^2$) for the group $\mathrm{E} = \Sigma(36 \times 3)$ is
\begin{equation}
\mathrm{d}\varepsilon = \frac{1}{48\pi^4} J^2 \, \mathrm{d}^{(4)} + \frac{1}{108\pi^4} J^2 \, \mathrm{d}^{(3)} + \frac{1}{3} \mathrm{d}^{(2)} - \frac{1}{18} \mathrm{d}^{(1)},
\end{equation}
where $\mathrm{d}^{(n)}$ is the uniform measure over the points in $D_n$.
\end{Thm}

\subsection{Group F $= \Sigma(72 \times 3)$}

The subgroup F has order 216, and its McKay graph, Figure \ref{fig:extended-E1(12)graph}, is the ``affine'' version of the graph $\mathcal{E}_1^{(12)}$ \cite[Figure 14]{evans/pugh:2009i}. It has a presentation with generators $S_1, S_2, T$ and $V$ of E, and the matrix $W \in SU(3)$ given by
$$W = \frac{1}{\sqrt{-3}} \left( \begin{array}{ccc} 1   &   1       & \omega^2  \\
                                                    1   & \omega    & \omega    \\
                                                 \omega &   1       & \omega    \end{array} \right).$$
The order of $W$ is 4. Note that $W^2 = V^2 S_1$. In fact, $S_1 = V^2 W^2$, $S_2 = (WV)^4$ and $T = V^2 W^3 V^2 W$ so that F $= \langle V, W \rangle$. It contains the group E as a normal subgroup.
This group is a subgroup of $SU(3)$ but not of $SU(3)/\mathbb{Z}_3$
The values of $\chi_{\rho}(\Gamma_j)$ for F are given in Table \ref{table:Character_table-(F)} (see \cite{hanany/he:1999}).

\renewcommand{\arraystretch}{1.5}

\begin{table}[bt]
\begin{center}
\small
\begin{tabular}{|c|ccc|ccc|} \hline
$j$ & 1 & 2 & 3 & 4, 5, 6 & 7, 8, 9 & 10, 11, 12 \\
\hline $|\Gamma_j|$ & 1 & 1 & 1 & 18 & 18 & 18 \\
\hline $\chi_{\rho}(\Gamma_j) \in \mathfrak{D}$ & 3 & $3\omega$ & $3\overline{\omega}$ & 1 & $\omega$ & $\overline{\omega}$ \\
\hline $(e^{2 \pi i \theta_1}, e^{2 \pi i \theta_2}) = \Phi^{-1}(\chi_{\rho}(\Gamma_j))$ & (1,1) & $(\omega, \overline{\omega})$ & $(\overline{\omega}, \omega)$ & $(1,i)$ & $(-\omega i,-\overline{\omega} i)$ & $(-\overline{\omega} i,-\omega i)$ \\
\hline $(\theta_1, \theta_2) \in [0,1]^2$ & (0,0) & $(\frac{1}{3}, \frac{2}{3})$ & $(\frac{2}{3}, \frac{1}{3})$ & $(0, \frac{1}{4})$ & $(\frac{1}{12}, \frac{5}{12})$ & $(\frac{5}{12}, \frac{1}{12})$ \\
\hline
\end{tabular}
$$\;$$
\begin{tabular}{|c|ccc|c|} \hline
$j$ &  13 & 14 & 15 & 16 \\
\hline $|\Gamma_j|$ & 9 & 9 & 9 & 24 \\
\hline $\chi_{\rho}(\Gamma_j) \in \mathfrak{D}$ & $-1$ & $-\omega$ & $-\overline{\omega}$ & 0 \\
\hline $(e^{2 \pi i \theta_1}, e^{2 \pi i \theta_2})  \in \mathbb{T}^2$ & (1,-1) & $(\omega,-\overline{\omega})$ & $(-\overline{\omega},\omega)$ & $(1,\omega)$ \\
\hline $(\theta_1, \theta_2) \in [0,1]^2$ & $(0, \frac{1}{2})$ & $(\frac{2}{6}, \frac{1}{6})$ & $(\frac{1}{6}, \frac{2}{6})$ & $(\frac{1}{3}, 0)$ \\
\hline
\end{tabular} \\
\normalsize
\caption{$\chi_{\rho}(\Gamma_j)$ for group F $= \Sigma(72 \times 3)$. Here $\omega = e^{2 \pi i/3}$.} \label{table:Character_table-(F)}
\end{center}
\end{table}

\renewcommand{\arraystretch}{1}

The spectral measure for the group F $= \Sigma(72 \times 3)$ is obtained in a similar way to that for the group E $= \Sigma(36 \times 3)$ in Section \ref{sect:(E)}, and we obtain:

\begin{Thm}
The spectral measure (over $\mathbb{T}^2$) for the group $\mathrm{F} = \Sigma(72 \times 3)$ is
\begin{equation}
\mathrm{d}\varepsilon = \frac{1}{32\pi^4} J^2 \, \mathrm{d}^{(4)} + \frac{1}{216\pi^4} J^2 \, \mathrm{d}^{(3)} + \frac{1}{6} \mathrm{d}^{(2)} - \frac{1}{36} \mathrm{d}^{(1)},
\end{equation}
where $\mathrm{d}^{(n)}$ is the uniform measure over the points in $D_n$.
\end{Thm}

\subsection{Group G $= \Sigma(216 \times 3)$}

The subgroup G has order 648, and its McKay graph, Figure \ref{fig:extended-E2(12)graph}, is the ``affine'' version of the graph $\mathcal{E}_2^{(12)}$ \cite[Figure 14]{evans/pugh:2009i}. It has a presentation with generators $S_1, S_2, T$ and $V$ of E, and the matrix $W \in SU(3)$ given by
$$U = \left( \begin{array}{ccc} \varepsilon^2   &   0           &   0   \\
                                    0           & \varepsilon^2 &   0   \\
                                    0           &   0           & \varepsilon^5  \end{array} \right),$$
where $\varepsilon = e^{2 \pi i/9}$.
The order of $U$ is 9. Note that $U^3 = S_2^2$.
It contains the group F as a normal subgroup.
This group is a subgroup of $SU(3)$ but not of $SU(3)/\mathbb{Z}_3$
The values of $\chi_{\rho}(\Gamma_j)$ for G are given in Table \ref{table:Character_table-(G)} (see \cite{desmier/sharp/patera:1982}).

\renewcommand{\arraystretch}{1.5}

\begin{table}[bt]
\begin{center}
\small
\begin{tabular}{|c|ccc|ccc|ccc|} \hline
$j$ & 1 & 2 & 3 & 4 & 5 & 6 & 7 & 8 & 9 \\
\hline $|\Gamma_j|$ & 1 & 1 & 1 & 54 & 54 & 54 & 9 & 9 & 9 \\
\hline $\chi_{\rho}(\Gamma_j) \in \mathfrak{D}$ & 3 & $3\omega$ & $3\overline{\omega}$ & 1 & $\omega$ & $\overline{\omega}$ & $-1$ & $-\omega$ & $-\overline{\omega}$ \\
\hline $(\theta_1, \theta_2) \in [0,1]^2$ & (0,0) & $(\frac{1}{3}, \frac{2}{3})$ & $(\frac{2}{3}, \frac{1}{3})$ & $(0, \frac{1}{4})$ & $(\frac{1}{12}, \frac{5}{12})$ & $(\frac{5}{12}, \frac{1}{12})$ & $(0, \frac{1}{2})$ & $(\frac{2}{6}, \frac{1}{6})$ & $(\frac{1}{6}, \frac{2}{6})$ \\
\hline
\end{tabular}
$$\;$$
\begin{tabular}{|c|ccc|ccc|} \hline
$j$ & 10 & 11 & 12 & 13 & 14 & 15 \\
\hline $|\Gamma_j|$ & 12 & 12 & 12 & 12 & 12 & 12 \\
\hline $\chi_{\rho}(\Gamma_j)$ & $\sqrt{3} e^{\pi i /18}$ & $\sqrt{3} e^{13 \pi i /18}$ & $\sqrt{3} e^{25 \pi i /18}$ & $\sqrt{3} e^{11 \pi i /18}$ & $\sqrt{3} e^{23 \pi i /18}$ & $\sqrt{3} e^{35 \pi i /18}$ \\
\hline $(\theta_1, \theta_2)$ & $(\frac{1}{9}, \frac{2}{9})$ & $(\frac{1}{9}, \frac{5}{9})$ & $(\frac{7}{9}, \frac{2}{9})$ & $(\frac{2}{9}, \frac{7}{9})$ & $(\frac{5}{9}, \frac{1}{9})$ & $(\frac{2}{9}, \frac{1}{9})$ \\
\hline
\end{tabular}
$$\;$$
\begin{tabular}{|c|ccc|ccc|c|c|} \hline
$j$ & 16 & 17 & 18 & 19 & 20 & 21 & 22 & 23, 24 \\
\hline $|\Gamma_j|$ & 36 & 36 & 36 & 36 & 36 & 36 & 24 & 72 \\
\hline $\chi_{\rho}(\Gamma_j)$ & $e^{5 \pi i/9}$ & $e^{11 \pi i/9}$ & $e^{17 \pi i/9}$ & $e^{\pi i/9}$ & $e^{7 \pi i/9}$ & $e^{13 \pi i/9}$ & 0 & 0 \\
\hline $(\theta_1, \theta_2)$ & $(\frac{2}{18}, \frac{7}{18})$ & $(\frac{8}{18}, \frac{1}{18})$ & $(\frac{5}{18}, \frac{1}{18})$ & $(\frac{1}{18}, \frac{5}{18})$ & $(\frac{1}{18}, \frac{8}{18})$ & $(\frac{7}{18}, \frac{2}{18})$ & $(0, \frac{1}{3})$ & $(0, \frac{1}{3})$ \\
\hline
\end{tabular} \\
\normalsize
\caption{$\chi_{\rho}(\Gamma_j)$ for group G $= \Sigma(216 \times 3)$. Here $\omega = e^{2 \pi i/3}$.} \label{table:Character_table-(G)}
\end{center}
\end{table}

\renewcommand{\arraystretch}{1}

The measures living on the points $(e^{2 \pi i \theta_1}, e^{2 \pi i \theta_2}) \in \mathbb{T}^2$ for $j=1,\ldots,9$ and $j=22,23,24$ have all been computed when we considered the subgroups E and F of $SU(3)$.

Let us denote be $\Sigma_{j_1}^{j_2}$ the summation
$$\Sigma_{j_1}^{j_2} = \frac{1}{6} \sum_{j=j_1}^{j_2} \sum_{g \in S_3} (g(\widetilde{\omega}^{\theta_1} + \widetilde{\omega}^{-\theta_2} + \widetilde{\omega}^{\theta_2-\theta_1}))^m (g(\widetilde{\omega}^{-\theta_1} + \widetilde{\omega}^{\theta_2} + \widetilde{\omega}^{\theta_1-\theta_2}))^n,$$
where for each $j$, $\theta_1, \theta_2$ are given in Table \ref{table:Character_table-(G)}, $\widetilde{\omega} = e^{2 \pi i/18}$, and the action of $g \in S_3$ on $(\widetilde{\omega}^{\theta_1} + \widetilde{\omega}^{-\theta_2} + \widetilde{\omega}^{\theta_2-\theta_1})$ is defined as follows: suppose $g(\widetilde{\omega}^{\theta_1},\widetilde{\omega}^{\theta_2}) = (\widetilde{\omega}^p,\widetilde{\omega}^q)$, then $g(\widetilde{\omega}^{\theta_1} + \widetilde{\omega}^{-\theta_2} + \widetilde{\omega}^{\theta_2-\theta_1}) = (\widetilde{\omega}^p + \widetilde{\omega}^{-q} + \widetilde{\omega}^{q-p})$.

For $j = 10,\ldots,15$, the pairs $(\theta_1,\theta_2) \in [0,1]^2$ are $(1/9,2/9)$, $(1/9,5/9)$, $(7/9,2/9)$, $(2/9,1/9)$, $(2/9,7/9)$, $(5/9,1/9)$. The $S_3$-orbits for these pairs contain three points each. These points all lie on the boundary of the fundamental domain $C$, and are obtained by taking the points $(\theta_1,\theta_2)$ such that $(e^{2 \pi i \theta_1}, e^{2 \pi i \theta_2}) \in D_{3/2}$ and then removing the points $(\theta_1,\theta_2)$ such that $\theta_k \in \{ 0,1/3,2/3 \}$, $k=1,2$.
Then we obtain $3\Sigma_{10}^{15} = 27 \int_{\mathbb{T}^2} R_{m_1,m_2}(\omega_1,\omega_2) \mathrm{d}^{(3)}(\omega_1,\omega_2) - 9 \int_{\mathbb{T}^2} R_{m_1,m_2}(\omega_1,\omega_2) \mathrm{d}_{3} \omega_1 \; \mathrm{d}_{3/2} \omega_2$.

\begin{figure}[bt]
\begin{center}
  \includegraphics[width=150mm]{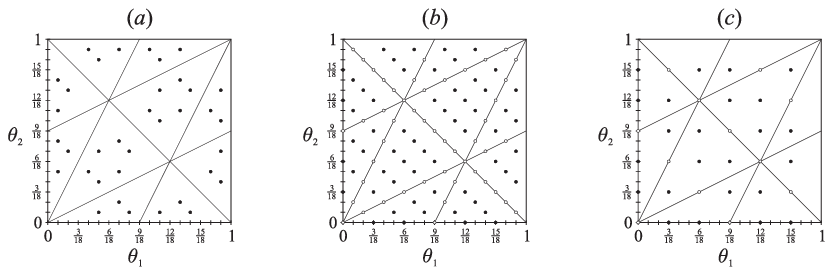}\\
 \caption{$(a)$ the points $\{ g(\theta_1,\theta_2) | g \in S_3 \}$ for $j = 16,\ldots,21$; $(b)$ the points $(\theta_1,\theta_2)$ such that $(e^{2 \pi i \theta_1}, e^{2 \pi i \theta_2}) \in D_6$; $(c)$ the points $(\theta_1,\theta_2)$ such that $e^{2 \pi i \theta_k}$ is a $6^{\mathrm{th}}$ root of unity, $k=1,2$.} \label{fig:poly-10}
\end{center}
\end{figure}

For $j=16,\ldots,21$, the pairs $(\theta_1,\theta_2)$ are $(5/18,1/18)$, $(2/18,7/18)$, $(8/18,1/18)$, $(1/18,5/18)$, $(1/18,8/18)$, $(7/18,2/18)$. The $S_3$-orbits for these pairs contain six points each, which are illustrated in Figure \ref{fig:poly-10}$(a)$. At these points, $J^2 = 48 \pi^4$. We can obtain this distribution by taking the points $(\theta_1,\theta_2)$ such that $(e^{2 \pi i \theta_1}, e^{2 \pi i \theta_2}) \in D_6$, illustrated in Figure \ref{fig:poly-10}$(b)$, each with the weight $J^2$ evaluated at that point. Since the points indicated by white circles in Figure \ref{fig:poly-10}$(b)$ map to the deltoid, the boundary of the discoid $\mathfrak{D}$, here $J^2 = 0$. We must then remove the points indicated by black circles in the interior of the triangular regions in Figure \ref{fig:poly-10}$(b)$ which are not in $\{ g(\theta_1,\theta_2) | g \in S_3 \}$.
Then we have $\Sigma_{16}^{21} = 108 \int_{\mathbb{T}^2} R_{m_1,m_2}(\omega_1,\omega_2) J^2 \, \mathrm{d}^{(6)}(\omega_1,\omega_2)/288\pi^4 - 36 \int_{\mathbb{T}^2} R_{m_1,m_2}(\omega_1,\omega_2) J^2 \, \mathrm{d}_3\omega_1 \; \mathrm{d}_3\omega_2/72\pi^4$.

Thus we obtain the following result:

\begin{Thm}
The spectral measure (over $\mathbb{T}^2$) for the group $\mathrm{G} = \Sigma(216 \times 3)$, is
\begin{eqnarray}
\mathrm{d}\varepsilon & = & \frac{1}{48\pi^4} J^2 \, \mathrm{d}^{(6)} + \frac{1}{96\pi^4} J^2 \, \mathrm{d}^{(4)} + \left( \frac{1}{6} + \frac{7}{648\pi^4} J^2 \right) \mathrm{d}^{(3)} + \frac{1}{18} \mathrm{d}^{(2)} - \frac{1}{108} \mathrm{d}^{(1)} \nonumber \\
& & - \frac{1}{36\pi^4} J^2 \, \mathrm{d}_3 \; \mathrm{d}_3 - \frac{1}{18} \mathrm{d}_{3/2} \; \mathrm{d}_{3/2},
\end{eqnarray}
where $\mathrm{d}_m$ is the uniform measure over $2m^{\mathrm{th}}$ roots of unity and $\mathrm{d}^{(m)}$ is the uniform measure on the points in $D_m$.
\end{Thm}

\subsection{Group H $\cong A_5$}

The subgroup H is the alternating $A_5$ group, which has order 60. Its McKay graph, Figure \ref{fig:extended-Hgraph}, is not the ``affine'' version of any of the $SU(3)$ $\mathcal{ADE}$ graphs. The values of $\chi_{\rho}(\Gamma_j)$  for H are given in Table \ref{table:Character_table-(H)}.

\renewcommand{\arraystretch}{1.5}

\begin{table}[bt]
\begin{center}
\small
\begin{tabular}{|c|c|c|c|cc|} \hline
$j$ & 1 & 2 & 3 & 4 & 5 \\
\hline $|\Gamma_j|$ & 1 & 20 & 15 & 12 & 12 \\
\hline $\chi_{\rho}(\Gamma_j) \in \mathfrak{D}$ & 3 & 0 & $-1$ & $\mu^+$ & $\mu^-$ \\
\hline $(e^{2 \pi i \theta_1}, e^{2 \pi i \theta_2}) = \Phi^{-1}(\chi_{\rho}(\Gamma_j)) \in \mathbb{T}^2$ & $(1,1)$ & $(1,\omega)$ & $(1,-1)$ & $(1,e^{\frac{2 \pi i}{5}})$ & $(1,e^{\frac{4 \pi i}{5}})$ \\
\hline $(\theta_1, \theta_2) \in [0,1]^2$ & (0,0) & $(\frac{1}{3}, \frac{1}{3})$ & $(0, \frac{1}{2})$ & $(0, \frac{1}{5})$ & $(0, \frac{2}{5})$ \\
\hline
\end{tabular} \\
\normalsize
\caption{$\chi_{\rho}(\Gamma_j)$ for group H $\cong A_5$. Here $\mu^{\pm} = (1 \pm \sqrt{5})/2$.} \label{table:Character_table-(H)}
\end{center}
\end{table}

\renewcommand{\arraystretch}{1}

For $j=1,3$, the points $(e^{2 \pi i \theta_1}, e^{2 \pi i \theta_2}) \in \mathbb{T}^2$ give the measure $(1/5)\mathrm{d}_1 \times \mathrm{d}_1 - (7/30)\delta_{(0,0)}$, where $\delta_x$ is the Dirac measure at the point $x$, whilst for $j=2$ we obtain the measure $J^2 \mathrm{d}^{(3)}/J^2\pi^4$ as in Section \ref{sect:Delta(3n^2)}. Thus we have the following result:

\begin{Thm}
The spectral measure (over $\mathbb{T}^2$) for the group $\mathrm{H} \cong A_5$, is
\begin{equation}
\mathrm{d}\varepsilon = \frac{1}{72\pi^4} J^2 \, \mathrm{d}^{(3)} + \frac{1}{5} \mathrm{d}_1 \times \mathrm{d}_1 - \frac{7}{30} \delta_{(0,0)} + \frac{1}{30} \sum_{k,l} \left( \delta_{(e^{2 \pi i k/5},e^{2 \pi i l/5})} + \delta_{(e^{4 \pi i k/5},e^{4 \pi i l/5})} \right),
\end{equation}
where the summation is over $k,l$ such that $(\pm k,\pm l) \in \{ (0,1), (1,0), (1,1) \}$.
The measure $\mathrm{d}_1$ is the uniform measure over 1 and $-1$, $\mathrm{d}^{(3)}$ is the uniform measure on the points in $D_3$ and $\delta_x$ is the Dirac measure at the point $x$.
\end{Thm}

\subsection{Group I $= \Sigma(168) \cong PSL(2,7)$}

The subgroup H is the projective special linear group $PSL(2,7)$, which has order 60. The special linear group $SL(2,7)$ consists of all $2 \times 2$ matrices over $\mathbb{F}_7$, the finite field with 7 elements, with unit determinant. The projective special linear group $PSL(2,7)$ is the quotient group $SL(2,7) / \{I,-I\}$, obtained by identifying $I$ and $-I$, where $I$ is the identity matrix. It is the automorphism group of the Klein quartic as well as the symmetry group of the Fano plane. Its McKay graph, Figure \ref{fig:extended-Igraph}, is not the ``affine'' version of any of the $SU(3)$ $\mathcal{ADE}$ graphs. The values of $\chi_{\rho}(\Gamma_j)$  for H are given in Table \ref{table:Character_table-(I)}.

\renewcommand{\arraystretch}{1.5}

\begin{table}[bt]
\begin{center}
\small
\begin{tabular}{|c|c|c|c|c|cc|} \hline
$j$ & 1 & 2 & 3 & 4 & 5 & 6 \\
\hline $|\Gamma_j|$ & 1 & 21 & 42 & 56 & 24 & 24 \\
\hline $\chi_{\rho}(\Gamma_j) \in \mathfrak{D}$ & 3 & $-1$ & 1 & 0 & $\nu$ & $\overline{\nu}$ \\
\hline $(e^{2 \pi i \theta_1}, e^{2 \pi i \theta_2}) = \Phi^{-1}(\chi_{\rho}(\Gamma_j))$ & $(1,1)$ & $(1,-1)$ & $(1,i)$ & $(1,\omega)$ & $(e^{\frac{2 \pi i}{7}},e^{\frac{6 \pi i}{7}})$ & $(e^{\frac{6 \pi i}{7}},e^{\frac{2 \pi i}{7}})$ \\
\hline $(\theta_1, \theta_2) \in [0,1]^2$ & (0,0) & $(0, \frac{1}{2})$ & $(0, \frac{1}{4})$ & $(\frac{1}{3}, \frac{1}{3})$ & $(\frac{1}{7}, \frac{3}{7})$ & $(\frac{3}{7}, \frac{1}{7})$ \\
\hline
\end{tabular} \\
\normalsize
\caption{$\chi_{\rho}(\Gamma_j)$ for group I $= \Sigma(168) \cong PSL(2,7)$. Here $\nu = (-1 + \sqrt{7}i)/2$.} \label{table:Character_table-(I)}
\end{center}
\end{table}

\renewcommand{\arraystretch}{1}

For $j=1,2$, the points $(e^{2 \pi i \theta_1}, e^{2 \pi i \theta_2}) \in \mathbb{T}^2$ give the measure $(1/6)\mathrm{d}_1 \times \mathrm{d}_1 - (1/28)\delta_{(0,0)}$, whilst for $j=3$ we obtain the measure $J^2 \mathrm{d}_2 \times \mathrm{d}_2$. For $j=4$ we obtain the measure $J^2 \mathrm{d}^{(3)}/J^2\pi^4$ as in Section \ref{sect:Delta(3n^2)}. Thus we have the following result:

\begin{Thm}
The spectral measure (over $\mathbb{T}^2$) for the group $\mathrm{I} = \Sigma(168) \cong PSL(2,7)$, is
\begin{equation}
\mathrm{d}\varepsilon = \frac{1}{72\pi^4} J^2 \, \mathrm{d}^{(3)} + \frac{1}{6} \mathrm{d}_1 \times \mathrm{d}_1 - \frac{1}{28} \delta_{(0,0)} + \frac{1}{42} \sum_{k,l} \left( \delta_{(e^{2 \pi i k/7},e^{2 \pi i l/7})} + \delta_{(e^{-2 \pi i k/7},e^{-2 \pi i l/7})} \right),
\end{equation}
where the summation is over $k$, $l$ such that $(k,l) \in \{ (1,3), (1,5), (2,3), (2,6), (4,5), (4,6) \}$.
The measure $\mathrm{d}_1$ is the uniform measure over 1 and $-1$, $\mathrm{d}^{(3)}$ is the uniform measure on the points in $D_3$ and $\delta_x$ is the Dirac measure at the point $x$.
\end{Thm}

\subsection{Group J $\cong TA_5$}

The subgroup J is the ternary alternating $A_5$ group, which has order 180. Its McKay graph, Figure \ref{fig:extended-Jgraph}, is not the ``affine'' version of any of the $SU(3)$ $\mathcal{ADE}$ graphs. The values of $\chi_{\rho}(\Gamma_j)$  for J are given in Table \ref{table:Character_table-(J)}. They are obtained from the character table of H $= \Sigma(60)$, Table \ref{table:Character_table-(H)}, where each conjugacy class $\Gamma'$ of the group H gives three conjugacy classes $\Gamma_j'$ of the group J, $j=1,2,3$, and $\chi_{\rho}(\Gamma_j') = \omega^{j-1} \lambda$ if $\chi_{\rho}(\Gamma') = \lambda$. Here $\omega = e^{2 \pi i/3}$ as usual.

\renewcommand{\arraystretch}{1.5}

\begin{table}[bt]
\begin{center}
\small
\begin{tabular}{|c|ccc|ccc|} \hline
$j$ & 1 & 2 & 3 & 4 & 5 & 6 \\
\hline $|\Gamma_j|$ & 1 & 1 & 1 & 15 & 15 & 15 \\
\hline $\chi_{\rho}(\Gamma_j) \in \mathfrak{D}$ & 3 & $3\omega$ & $3\overline{\omega}$ & $-1$ & $-\omega$ & $-\overline{\omega}$ \\
\hline $(e^{2 \pi i \theta_1}, e^{2 \pi i \theta_2}) = \Phi^{-1}(\chi_{\rho}(\Gamma_j)) \in \mathbb{T}^2$ & $(1,1)$ & $(\omega, \overline{\omega})$ & $(\overline{\omega}, \omega)$ & $(1,-1)$ & $(\omega,-\overline{\omega})$ & $(-\overline{\omega},\omega)$ \\
\hline $(\theta_1, \theta_2) \in [0,1]^2$ & (0,0) & $(\frac{1}{3}, \frac{2}{3})$ & $(\frac{2}{3}, \frac{1}{3})$ & $(0, \frac{1}{2})$ & $(\frac{2}{6}, \frac{1}{6})$ & $(\frac{1}{6}, \frac{2}{6})$ \\
\hline
\end{tabular}
$$\;$$
\begin{tabular}{|c|ccc|ccc|c|} \hline
$j$ & 7 & 8 & 9 & 10 & 11 & 12 & 13, 14, 15 \\
\hline $|\Gamma_j|$ & 12 & 12 & 12 & 12 & 12 & 12 & 20 \\
\hline $\chi_{\rho}(\Gamma_j)$ & $\mu^+$ & $\mu^+ \omega$ & $\mu^+ \overline{\omega}$ & $\mu^-$ & $\mu^- \omega$ & $\mu^- \overline{\omega}$ & 0 \\
\hline $(e^{2 \pi i \theta_1}, e^{2 \pi i \theta_2})$ & $(1,e^{\frac{2 \pi i}{5}})$ & $(\omega,e^{\frac{14 \pi i}{15}})$ & $(e^{\frac{14 \pi i}{15}},\omega)$ & $(1,e^{\frac{4 \pi i}{5}})$ & $(\omega,e^{\frac{8 \pi i}{15}})$ & $(e^{\frac{8 \pi i}{15}},\omega)$ & $(1,\omega)$ \\
\hline $(\theta_1, \theta_2)$ & $(0, \frac{1}{5})$ & $(\frac{2}{3}, \frac{7}{15})$ & $(\frac{7}{15}, \frac{2}{3})$ & $(0, \frac{2}{5})$ & $(\frac{2}{3}, \frac{4}{15})$ & $(\frac{4}{15}, \frac{2}{3})$ & $(\frac{1}{3}, 0)$ \\
\hline
\end{tabular} \\
\normalsize
\caption{$\chi_{\rho}(\Gamma_j)$ for group J $\cong TA_5$. Here $\omega = e^{2 \pi i/3}$ and $\mu^{\pm} = (1 \pm \sqrt{5})/2$.} \label{table:Character_table-(J)}
\end{center}
\end{table}

\renewcommand{\arraystretch}{1}

For $j=7,8,9$, the points $(e^{2 \pi i \theta_1}, e^{2 \pi i \theta_2}) \in \mathbb{T}^2$ give the measure $\mathrm{d}^{((5))}$, whilst for $j=10,11,12$ they give the measure $\mathrm{d}^{((5/2))}$, each with weight $1/5$. Thus we obtain the following result:

\begin{Thm}
The spectral measure (over $\mathbb{T}^2$) for the group $\mathrm{J} \cong TA_5$, is
\begin{equation}
\mathrm{d}\varepsilon = \frac{1}{72\pi^4} J^2 \, \mathrm{d}^{(3)} + \frac{1}{3} \mathrm{d}^{(2)} - \frac{1}{15} \mathrm{d}^{(1)} + \frac{1}{5} \mathrm{d}^{((5))} + \frac{1}{5} \mathrm{d}^{((5/2))},
\end{equation}
where the measures $\mathrm{d}_n$, $\mathrm{d}^{(n)}$, $\mathrm{d}^{((n))}$ are as in Definition \ref{def:4measures}.
\end{Thm}

\subsection{Group K $\cong TPSL(2,7)$}

The subgroup K is the ternary $PSL(2,7)$ group, and has order 504. Its McKay graph, Figure \ref{fig:extended-Kgraph}, is the ``affine'' version of the graph $\mathcal{E}_5^{(12)}$ \cite[Figure 15]{evans/pugh:2009i}. The values of $\chi_{\rho}(\Gamma_j)$ for K are given in Table \ref{table:Character_table-(K)}. They are obtained from the character table of I $= \Sigma(168)$, Table \ref{table:Character_table-(I)}, in the same way as the values of $\chi_{\rho}(\Gamma_j)$ for J $\cong TA_5$ are obtained from the character table of H $= \Sigma(60)$.

\renewcommand{\arraystretch}{1.5}

\begin{table}[bt]
\begin{center}
\small
\begin{tabular}{|c|ccc|ccc|} \hline
$j$ & 1 & 2 & 3 & 4 & 5 & 6 \\
\hline $|\Gamma_j|$ & 1 & 1 & 1 & 21 & 21 & 21 \\
\hline $\chi_{\rho}(\Gamma_j) \in \mathfrak{D}$ & 3 & $3\omega$ & $3\overline{\omega}$ & $1$ & $\omega$ & $\overline{\omega}$ \\
\hline $(e^{2 \pi i \theta_1}, e^{2 \pi i \theta_2}) = \Phi^{-1}(\chi_{\rho}(\Gamma_j))$ & $(1,1)$ & $(\omega, \overline{\omega})$ & $(\overline{\omega}, \omega)$ & $(1,i)$ & $(-\omega i,-\overline{\omega} i)$ & $(-\overline{\omega} i,-\omega i)$ \\
\hline $(\theta_1, \theta_2) \in [0,1]^2$ & (0,0) & $(\frac{1}{3}, \frac{2}{3})$ & $(\frac{2}{3}, \frac{1}{3})$ & $(0, \frac{1}{4})$ & $(\frac{1}{12}, \frac{5}{12})$ & $(\frac{5}{12}, \frac{1}{12})$ \\
\hline
\end{tabular}
$$\;$$
\begin{tabular}{|c|ccc|ccc|} \hline
$j$ & 7 & 8 & 9 & 10 & 11 & 12 \\
\hline $|\Gamma_j|$ & 42 & 42 & 42 & 24 & 24 & 24 \\
\hline $\chi_{\rho}(\Gamma_j)$ & -1 & $-\omega$ & $-\overline{\omega}$ & $\nu$ & $\nu \omega$ & $\nu \overline{\omega}$ \\
\hline $(e^{2 \pi i \theta_1}, e^{2 \pi i \theta_2})$ & $(1,-1)$ & $(\omega,-\overline{\omega})$ & $(-\overline{\omega},\omega)$ & $(e^{\frac{2 \pi i}{7}},e^{\frac{6 \pi i}{7}})$ & $(e^{\frac{20 \pi i}{21}},e^{\frac{16 \pi i}{21}})$ & $(e^{\frac{10 \pi i}{21}},e^{\frac{8 \pi i}{21}})$ \\
\hline $(\theta_1, \theta_2)$ & $(0, \frac{1}{2})$ & $(\frac{2}{6}, \frac{1}{6})$ & $(\frac{1}{6}, \frac{2}{6})$ & $(\frac{1}{7}, \frac{3}{7})$ & $(\frac{10}{21}, \frac{8}{21})$ & $(\frac{5}{21}, \frac{4}{21})$ \\
\hline
\end{tabular}
$$\;$$
\begin{tabular}{|c|ccc|c|} \hline
$j$ & 13 & 14 & 15 & 16, 17, 18 \\
\hline $|\Gamma_j|$ & 24 & 24 & 24 & 56 \\
\hline $\chi_{\rho}(\Gamma_j)$ & $\overline{\nu}$ & $\overline{\nu} \omega$ & $\overline{\nu \omega}$ & 0 \\
\hline $(e^{2 \pi i \theta_1}, e^{2 \pi i \theta_2})$ & $(e^{\frac{6 \pi i}{7}},e^{\frac{2 \pi i}{7}})$ & $(e^{\frac{8 \pi i}{21}},e^{\frac{10 \pi i}{21}})$ & $(e^{\frac{16 \pi i}{21}},e^{\frac{20 \pi i}{21}})$ & $(1,\omega)$ \\
\hline $(\theta_1, \theta_2)$ & $(\frac{3}{7}, \frac{1}{7})$ & $(\frac{4}{21}, \frac{5}{21})$ & $(\frac{8}{21}, \frac{10}{21})$ & $(\frac{1}{3}, 0)$ \\
\hline
\end{tabular} \\
\normalsize
\caption{$\chi_{\rho}(\Gamma_j)$ for group K $\cong TPSL(2,7)$. Here $\omega = e^{2 \pi i/3}$ and $\nu = (-1 + \sqrt{7}i)/2$.} \label{table:Character_table-(K)}
\end{center}
\end{table}

\renewcommand{\arraystretch}{1}

The measures which give the points $(e^{2 \pi i \theta_1}, e^{2 \pi i \theta_2}) \in \mathbb{T}^2$ for $j=1,\ldots,9$ and $j=16,17,18$ have all been computed for previous subgroups of $SU(3)$. For $j=10,\ldots,15$ we obtain the measure $\mathrm{d}^{(1/21,21/4)}$. Thus we have:

\begin{Thm}
The spectral measure (over $\mathbb{T}^2$) for the group $\mathrm{K} \cong TPSL(2,7)$, is
\begin{equation}
\mathrm{d}\varepsilon = \frac{1}{32\pi^4} J^2 \, \mathrm{d}^{(4)} + \frac{1}{24\pi^4} J^2 \, \mathrm{d}^{(3)} + \frac{1}{2} \mathrm{d}^{(2)} - \frac{5}{42} \mathrm{d}^{(1)} + \frac{6}{7} \mathrm{d}^{(1/21,21/4)},
\end{equation}
where the measures $\mathrm{d}_n$, $\mathrm{d}^{(n)}$, $\mathrm{d}^{(n,k)}$ are as in Definition \ref{def:4measures}.
\end{Thm}

\subsection{Group L $= \Sigma(360 \times 3) \cong TA_6$}

\renewcommand{\arraystretch}{1.5}

\begin{table}[bt]
\begin{center}
\small
\begin{tabular}{|c|ccc|ccc|} \hline
$j$ & 1 & 2 & 3 & 4 & 5 & 6 \\
\hline $|\Gamma_j|$ & 1 & 1 & 1 & 90 & 90 & 90 \\
\hline $\chi_{\rho}(\Gamma_j) \in \mathfrak{D}$ & 3 & $3\omega$ & $3\overline{\omega}$ & $1$ & $\omega$ & $\overline{\omega}$ \\
\hline $(e^{2 \pi i \theta_1}, e^{2 \pi i \theta_2}) = \Phi^{-1}(\chi_{\rho}(\Gamma_j))$ & $(1,1)$ & $(\omega, \overline{\omega})$ & $(\overline{\omega}, \omega)$ & $(1,i)$ & $(-\omega i,-\overline{\omega} i)$ & $(-\overline{\omega} i,-\omega i)$ \\
\hline $(\theta_1, \theta_2) \in [0,1]^2$ & (0,0) & $(\frac{1}{3}, \frac{2}{3})$ & $(\frac{2}{3}, \frac{1}{3})$ & $(0, \frac{1}{4})$ & $(\frac{1}{12}, \frac{5}{12})$ & $(\frac{5}{12}, \frac{1}{12})$ \\
\hline
\end{tabular}
$$\;$$
\begin{tabular}{|c|ccc|ccc|} \hline
$j$ & 7 & 8 & 9 & 10 & 11 & 12 \\
\hline $|\Gamma_j|$ & 45 & 45 & 45 & 72 & 72 & 72 \\
\hline $\chi_{\rho}(\Gamma_j)$ & -1 & $-\omega$ & $-\overline{\omega}$ & $\mu^+$ & $\mu^+ \omega$ & $\mu^+ \overline{\omega}$ \\
\hline $(e^{2 \pi i \theta_1}, e^{2 \pi i \theta_2})$ & $(1,-1)$ & $(\omega,-\overline{\omega})$ & $(-\overline{\omega},\omega)$ & $(1,e^{\frac{2 \pi i}{5}})$ & $(\omega,e^{\frac{14 \pi i}{15}})$ & $(e^{\frac{14 \pi i}{15}},\omega)$ \\
\hline $(\theta_1, \theta_2)$ & $(0, \frac{1}{2})$ & $(\frac{2}{6}, \frac{1}{6})$ & $(\frac{1}{6}, \frac{2}{6})$ & $(0, \frac{1}{5})$ & $(\frac{2}{3}, \frac{7}{15})$ & $(\frac{7}{15}, \frac{2}{3})$ \\
\hline
\end{tabular}
$$\;$$
\begin{tabular}{|c|ccc|c|} \hline
$j$ & 13 & 14 & 15 & 16, 17 \\
\hline $|\Gamma_j|$ & 72 & 72 & 72 & 120 \\
\hline $\chi_{\rho}(\Gamma_j) \in \mathfrak{D}$ & $\mu^-$ & $\mu^- \omega$ & $\mu^- \overline{\omega}$ & 0 \\
\hline $(e^{2 \pi i \theta_1}, e^{2 \pi i \theta_2}) \in \mathbb{T}^2$ & $(1,e^{\frac{4 \pi i}{5}})$ & $(\omega,e^{\frac{8 \pi i}{15}})$ & $(e^{\frac{8 \pi i}{15}},\omega)$ & $(1,\omega)$ \\
\hline $(\theta_1, \theta_2) \in [0,1]^2$ & $(0, \frac{2}{5})$ & $(\frac{2}{3}, \frac{4}{15})$ & $(\frac{4}{15}, \frac{2}{3})$ & $(\frac{1}{3}, 0)$ \\
\hline
\end{tabular} \\
\normalsize
\caption{$\chi_{\rho}(\Gamma_j)$ for group L $= \Sigma(360 \times 3) \cong TA_6$. Here $\omega = e^{2 \pi i/3}$ and $\mu^{\pm} = (1 \pm \sqrt{5})/2$.} \label{table:Character_table-(L)}
\end{center}
\end{table}

\renewcommand{\arraystretch}{1}

The subgroup L is the ternary alternating $A_6$ group, which has order 1080. Its McKay graph, Figure \ref{fig:extended-E4(12)graph}, is the ``affine'' version of the graph $\mathcal{E}_4^{(12)}$ \cite[Figure 7]{evans/pugh:2009ii}. The values of $\chi_{\rho}(\Gamma_j)$ for L are given in Table \ref{table:Character_table-(L)} (see \cite{desmier/sharp/patera:1982}).
The values of $\chi_{\rho}(\Gamma_j)$ for the group L $= \Sigma(360 \times 3) \cong TA_6$ have all appeared for previous groups, and hence it is easy to compute the spectral measure:

\begin{Thm}
The spectral measure (over $\mathbb{T}^2$) for the group $\mathrm{L} = \Sigma(360 \times 3) \cong TA_6$, is
\begin{equation}
\mathrm{d}\varepsilon = \frac{1}{96\pi^4} J^2 \, \mathrm{d}^{(4)} + \frac{1}{108\pi^4} J^2 \, \mathrm{d}^{(3)} + \frac{1}{6} \mathrm{d}^{(2)} - \frac{7}{180} \mathrm{d}^{(1)} + \frac{1}{5} \mathrm{d}^{((5))} + \frac{1}{5} \mathrm{d}^{((5/2))},
\end{equation}
where the measures $\mathrm{d}_n$, $\mathrm{d}^{(n)}$, $\mathrm{d}^{((n))}$ are as in Definition \ref{def:4measures}.
\end{Thm}

\paragraph{Acknowledgements}

This work was supported by the Marie Curie Research Training Network MRTN-CT-2006-031962 EU-NCG.

\end{document}